\documentclass[final,onefignum,onetabnum]{siamart190516}

\usepackage{lipsum}
\usepackage{amsfonts}
\usepackage{graphicx}
\usepackage{epstopdf}
\usepackage{algorithmic}
\usepackage{arydshln}
\usepackage{enumitem}

\ifpdf
  \DeclareGraphicsExtensions{.eps,.pdf,.png,.jpg}
\else
  \DeclareGraphicsExtensions{.eps}
\fi

\newcommand{\define}{\stackrel{\text{\tiny def}}{=}}

\newsiamremark{remark}{Remark}
\newsiamremark{hypothesis}{Hypothesis}
\crefname{hypothesis}{Hypothesis}{Hypotheses}
\newsiamthm{claim}{Claim}

\headers{Tensor Expander Chernoff Bounds}{Shih Yu Chang}

\title{Tensor Expander Chernoff Bounds}

\author{Shih Yu Chang\thanks{Department of Applied Data Science, San Jose State University, San Jose, CA 95192, USA 
  (\email{ shihyu.chang@sjsu.edu } ).}
}

\usepackage{amsopn}

\ifpdf
\hypersetup{
  pdftitle={Tensor Expander Chernoff Bounds},
  pdfauthor={Shih Yu Chang}
}
\fi

\begin{document}

\maketitle


\begin{abstract}
The Chernoff bound is an important inequality relation in probability theory. The original version of the Chernoff bound is to give an exponential decreasing bound on the tail distribution of sums of independent random variables. Recent years, several works have been done by extending the original version of the Chernoff bound to high-dimensional random objects, e.g., random matrices, or/and to consider the relaxation that there is no requirement of independent assumptions among random objects. In this work, we generalize the matrix expander Chernoff bound studied by Garg et al.~\cite{garg2018matrix} to \emph{tensor expander Chernoff bounds}. Our main tool is to develop new tensor norm inequalities based on log-majorization techniques. These new tensor norm inequalities are used to bound the expectation of Ky Fan norm of the random tensor exponential function, then \emph{tensor expander Chernoff bounds} can be established. Compared with the matrix expander Chernoff bound, the \emph{tensor expander Chernoff bounds} proved at this work contributes following aspects: (1) the random objects dimensions are increased from matrices (two-dimensional data array) to tensors (multidimensional data array); (2) this bound generalizes the identity map of the random objects summation to any polynomial function of the random objects summation; (3) Ky Fan norm, instead only the maximum or the minimum eigenvalues, for the function of the random objects summation is considered; (4) we remove the restriction about the summation of all mapped random objects is zero, which is required in the matrix expander Chernoff bound derivation. 
\end{abstract}



\section{Introduction}




In probability theory, the Chernoff bound provides the exponential decreasing bound on tail distribution of sums of independent random variables. This bound has various applications in science and engineering. For example, the Chernoff bound is utilized in mathematical learning theory to prove that a learning algorithm is probably approximately correct~\cite{kearns1994introduction}. On the engineering side, the Chernoff bound is also used to obtain tight bounds for packet routing problems which reduce wireless communication congestion while routing packets in sparse networks~\cite{jang1999chernoff}.

It is a tighter bound than the known first- or second-moment-based tail bounds such as Markov's inequality or Chebyshev's inequality, which only yield power-law bounds on tail distribution. Nevertheless, neither Markov's inequality nor Chebyshev's inequality requires that the variates are independent, which is necessary by the Chernoff bound~\cite{chernoff1981note}. Given $n$ independent and identically distributed (i.i.d.) random variables $X_1, X_2, \cdots, X_n$ taking values in $\{0, 1\}$ with $\mathbb{E}[X_i] = q$ and $\epsilon > 0$, Chernoff bound for the version of $n$ i.i.d. random variables is
\begin{eqnarray}\label{eq:Chernoff bound iid rvs}
\mathrm{Pr} \left(\frac{1}{n} \sum\limits_{i=1}^{n}X_i \geq q + \epsilon \right) \leq 
\exp (- n \mathbb{D} ( q + \epsilon \parallel q )),
\end{eqnarray}
where $\mathbb{D}(x \parallel y) \define x \ln \frac{x}{y} + (1 - x) \ln \frac{1-x}{1 - y}$ is the Kullback–Leibler divergence between Bernoulli distributed random variables with parameters $x$ and $y$ respectively.

There are various directions to generalize the Chernoff bound from Eq.~\eqref{eq:Chernoff bound iid rvs}, one major direction is to increase the dimension of random objects from random variables to random matrices. The works of Rudelson~\cite{rudelson1999random}, Ahlswede-Winter~\cite{ahlswede2002strong} and Tropp~\cite{tropp2012user} demonstrated that a similar concentration bound is also valid for matrix-valued random variables. If $\mathbf{X}_1,\mathbf{X}_2, \cdots, \mathbf{X}_n$ are independent $m \times m$ Hermitian complex random matrices with $\left\Vert  \mathbf{X}_i \right\Vert \leq 1$ for $1 \leq i \leq n$, where $\left\Vert \cdot \right\Vert$ is the spectral norm, we have following Chernoff bound for the version of $n$ i.i.d. random matrices: 
\begin{eqnarray}\label{eq:Chernoff bound iid matrices}
\mathrm{Pr} \left( \left\Vert \frac{1}{n} \sum\limits_{i=1}^{n} \mathbf{X}_i - \mathbb{E}[\mathbf{X}] \right\Vert \geq \vartheta \right) \leq m \exp( - \Omega  n \vartheta^2),
\end{eqnarray}
where $\Omega$ is a constant related to the matrix norm. This is also called ``Matrix Chernoff Bound'' and is applied to many fields, e.g., spectral graph theory, numerical linear algebra, machine learning and information theory~\cite{tropp2015introduction}. Recently, Shih Yu generalized matrix bounds to various tensors bounds, e.g., Chernoff, Bennett, and Bernstein inequalities associated with tensors in~\cite{chang2020convenient}. 

Another direction to extend from the basic Chernoff bound is to consider non-independent assumptions for random variables. By Gillman~\cite{gillman1998chernoff} and its refinement works~\cite{chung2012chernoff, rao2017sharp}, they changed the independence assumption to Markov dependence and we summarize their works as follows. We are given $\mathfrak{G}$ as a regular $\lambda$-expander graph with vertex set $\mathfrak{V}$, and $g: \mathfrak{V} \rightarrow \mathbb{C}$ as a bounded function. Suppose $v_1, v_2\cdots, v_{\kappa}$ is a stationary random walk of length $\kappa$ on $\mathfrak{G}$, it is shown that:
\begin{eqnarray}\label{eq:Chernoff bound Markov rvs}\mathrm{Pr} \left( \left\Vert \frac{1}{\kappa} \sum\limits_{j=1}^{\kappa} g(v_i) - \mathbb{E}[g] \right\Vert \geq \vartheta \right) \leq 2 \exp( - \Omega (1 - \lambda) \kappa \vartheta^2).
\end{eqnarray}
The value of $\lambda$ is also the second largest eigenvalue of the transition matrix of the underlying graph $\mathfrak{G}$. The bound given in Eq.~\eqref{eq:Chernoff bound Markov rvs} is named as ``Expander Chernoff Bound''. It is natural to generalize Eq.~\eqref{eq:Chernoff bound Markov rvs} to ``Matrix Expander Chernoff Bound''. Wigderson and Xiao in~\cite{wigderson2008derandomizing} began first attempt to obtain partial results of ``Matrix Expander Chernoff Bound'' and the complete solution is given later by Garg et al.~\cite{garg2018matrix}. Their results can be summarized by the following theorem.
\begin{theorem}\label{thm:matrix expander Chernoff bound}
Let $\mathfrak{G} = (\mathfrak{V}, \mathfrak{E})$ be a regular graph whose transition matrix has second largest eigenvalue as $\lambda$, and let $g: \mathfrak{V} \rightarrow \mathbb{C}^{m \times m}$ be a function satisfy following:
\begin{enumerate}
\item For each $v \mathfrak{V}$, $g(v)$ is a Hermitian matrix with $ \left\Vert g(v) \right\Vert \leq 1$;
\item $\sum\limits_{v \in \mathfrak{V}} g(v) = \mathbf{0}$.
\end{enumerate}
Then, for a stationary random walk $v_1, \cdots, v_\kappa$ with $\epsilon \in (0, 1)$, we have
\begin{eqnarray}\label{eq:thm:matrix expander Chernoff bound}
\mathrm{P}\left( \lambda_{g, \max} \left(\frac{1}{\kappa} \sum\limits_{j=1}^{\kappa} g(v_j)\right) \geq \epsilon\right) & \leq &  m \exp(- \Omega (1 - \lambda)\kappa \epsilon^2), \nonumber \\
\mathrm{P}\left( \lambda_{g, \min} \left(\frac{1}{\kappa} \sum\limits_{j=1}^{\kappa} g(v_j)\right) \leq - \epsilon\right) & \leq &  m \exp(- \Omega (1 - \lambda)\kappa \epsilon^2),
\end{eqnarray}
where $\lambda_{g, \max}, (\lambda_{g, \min})$ is the largest (smallest) eigenvalue of the summation of $\kappa$ matrices obtained by the mapping $g$. 
\end{theorem}


In this work, we generalize matrix expander Chernoff bound to \emph{tensor expander Chernoff bounds} by allowing more general norm for tensors, Ky Fan norm and general convex function, instead identity function, of the tensors summand. We first extend the setting of Hermitian matrices to Hermitian tensors in~\cite{hiai2017generalized} by utilizing majorization techniques to prove the following theorem~\ref{thm:Multivaraite Tensor Norm Inequalities intro}, which will play a crucial role in our proof of \emph{tensor expander Chernoff bounds}. 
\begin{theorem}\label{thm:Multivaraite Tensor Norm Inequalities intro}
Let $\mathcal{C}_i \in \mathbb{C}^{I_1 \times \cdots \times I_N \times I_1 \times \cdots \times I_N}$ be positive Hermitian tensors for $1 \leq i \leq n$ with Hermitian rank $r$, $\left\Vert \cdot \right\Vert_{(k)}$ be a Ky Fan $k$-norm with corresponding gauge function $\rho$. For any continuous function $f:(0, \infty) \rightarrow [0, \infty)$ such that $x \rightarrow \log f(e^x)$ is convex on $\mathbb{R}$, we have 
\begin{eqnarray}\label{eq1:thm:Multivaraite Tensor Norm Inequalities intro}
\left\Vert  f \left( \exp \left( \sum\limits_{i=1}^n \log \mathcal{C}_i\right)   \right)  \right\Vert_{(k)} &\leq& \exp \int_{- \infty}^{\infty} \log \left\Vert f \left( \left\vert \prod\limits_{i=1}^{n}  \mathcal{C}_i^{1 + \iota t} \right\vert\right)\right\Vert_{(k)} \beta_0(t) dt ,
\end{eqnarray}
where $\iota$ is $\sqrt{-1}$ and $\beta_0(t) = \frac{\pi}{2 (\cosh (\pi t) + 1)}$. For any continuous function $g(0, \infty) \rightarrow [0, \infty)$ such that $x \rightarrow g (e^x)$ is convex on $\mathbb{R}$, we have 
\begin{eqnarray}\label{eq2:thm:Multivaraite Tensor Norm Inequalities intro}
\left\Vert  g \left( \exp \left( \sum\limits_{i=1}^n \log \mathcal{C}_i\right)   \right)  \right\Vert_{(k)} &\leq& \int_{- \infty}^{\infty} \left\Vert g \left( \left\vert \prod\limits_{i=1}^{n}  \mathcal{C}_i^{1 + \iota t} \right\vert\right)\right\Vert_{(k)} \beta_0(t) dt. 
\end{eqnarray}
\end{theorem}

The main contribution of this work is the following \emph{tensor expander Chernoff bounds}. 
\begin{theorem}\label{thm:tensor expander intro}
Let $\mathfrak{G} = (\mathfrak{V}, \mathfrak{E})$ be a regular undirected graph whose transition matrix has second eigenvalue $\lambda$, and let $g: \mathfrak{V} \rightarrow \in \mathbb{C}^{I_1 \times \cdots \times I_M \times I_1 \times \cdots \times I_M}$ be a function. We assume following: 
\begin{enumerate}
\item For each $v \in \mathfrak{V}$, $g(v)$ is a Hermitian tensor;
\item $\left\Vert g(v) \right\Vert \leq r$;
\item A nonnegative coefficients polynomial raised by the power $s \geq 1$ as $f: x \rightarrow (a_0 + a_1x +a_2 x^2 + \cdots +a_n x^n)^s$ satisfying $f \left(\exp \left( t   \sum\limits_{j=1}^{\kappa} g(v_j) \right) \right) \succeq \exp \left( t f \left(  \sum\limits_{j=1}^{\kappa} g(v_j) \right) \right) $ almost surely;
\item  For $\tau \in [\infty, \infty]$, we have constants $C$ and $\sigma$ such that $ \beta_0(\tau)  \leq \frac{C \exp( \frac{-\tau^2}{2 \sigma^2} ) }{\sigma \sqrt{2 \pi}}$. 
\end{enumerate}
Then, we have 
\begin{eqnarray}\label{eq0:thm:tensor expander intro}
\mathrm{Pr}\left( \left\Vert f \left(  \sum\limits_{j=1}^{\kappa} g(v_j) \right)  \right\Vert_{(k)} \geq\vartheta \right) \leq   \min\limits_{t > 0 } \left[ (n+1)^{(s-1)} e^{-\vartheta t} \left(a_0 k  +C \left( k + \sqrt{\frac{\mathbb{I}_1^M - k }{k}}\right)\cdot \right. \right. ~~~~~  \nonumber \\
\left. \left. \sum\limits_{l=1}^n a_l \exp( 8 \kappa \overline{\lambda} + 2  (\kappa +8 \overline{\lambda}) lsr t + 2(\sigma (\kappa +8 \overline{\lambda}) lsr )^2 t^2  )  \right)\right],
\end{eqnarray}
where $\mathbb{I}_1^M$ is a positive integer obtained from $\prod\limits_{k=1}^M I_k$. 
\end{theorem}
By comparing matrix expander Chernoff bound, the \emph{tensor expander Chernoff bounds} derived at this paper makes following relaxation: (1) the random objects dimensions are increased from matrices (2D data array) to tensors (multidimensional data array); (2)this bound generalizes the identity map to the power of polynomial functions as shown by the third assumption of the function $f$ at Theorem~\ref{thm:tensor expander intro}; (3) Ky Fan norm of $ f \left(  \sum\limits_{j=1}^{\kappa} g(v_j) \right) $  is considered instead only the maximum or the minimum eigenvalues of $ f \left(  \sum\limits_{j=1}^{\kappa} g(v_j) \right) $ being evaluated; (4) there are no restriction for $\sum\limits_{v \in \mathfrak{V}} g(v)$, but this summation is required to be a zero matrix in the matrix expander Chernoff bound derivation~\footnote{We have another work~\cite{chang2021general} apply same majorization techniques to build bounds for tail bounds for the summation of random tensors under indepedent assumptions.}.


The paper is organized as follows. Preliminaries of tensors and basic majorization notations are given in Section~\ref{sec:Fundamentals of Tensors and Majorization}. In Section~\ref{sec:Multivariate Tensor Norm Inequalities}, multivariate tensor norm inequalities are established and these inequalities will become our main ingredients to prove \emph{tensor expander Chernoff bounds}. Main theorem~\ref{thm:tensor expander intro} is discussed and proved in Section~\ref{sec:Tensor Expander Chernoff Bounds Derivation by Majorization}. Finally, the conclusio and potential future works are given in Section~\cref{sec:Conclusions}. 

\section{Fundamentals of Tensors and Majorization}\label{sec:Fundamentals of Tensors and Majorization}

The purpose of this section is to provide fundamental facts about tensors and introduce notions about majorization. 

\subsection{Tensors Preliminaries}\label{sec:Tensors Preliminaries}

Throughout this work, scalars are represented by lower-case letters (e.g., $d$, $e$, $f$, $\ldots$), vectors by boldfaced lower-case letters (e.g., $\mathbf{d}$, $\mathbf{e}$, $\mathbf{f}$, $\ldots$), matrices by boldfaced capitalized letters (e.g., $\mathbf{D}$, $\mathbf{E}$, $\mathbf{F}$, $\ldots$), and tensors by calligraphic letters (e.g., $\mathcal{D}$, $\mathcal{E}$, $\mathcal{F}$, $\ldots$), respectively. Tensors are multiarrays of values which are higher-dimensional generalizations from vectors and matrices. Given a positive integer $N$, let $[N] \define \{1, 2, \cdots ,N\}$. An \emph{order-$N$ tensor} (or \emph{$N$-th order tensor}) denoted by $\mathcal{X} \define (x_{i_1, i_2, \cdots, i_N})$, where $1 \leq i_j = 1, 2, \ldots, I_j$ for $j \in [N]$, is a multidimensional array containing $\prod_{n=1}^N I_n$ entries. 
Let $\mathbb{C}^{I_1 \times \cdots \times I_N}$ and $\mathbb{R}^{I_1 \times \cdots \times I_N}$ be the sets of the order-$N$ $I_1 \times \cdots \times I_N$ tensors over the complex field $\mathbb{C}$ and the real field $\mathbb{R}$, respectively. For example, $\mathcal{X} \in \mathbb{C}^{I_1 \times \cdots \times I_N}$ is an order-$N$ multiarray, where the first, second, ..., and $N$-th dimensions have $I_1$, $I_2$, $\ldots$, and $I_N$ entries, respectively. Thus, each entry of $\mathcal{X}$ can be represented by $x_{i_1, \cdots, i_N}$. For example, when $N = 3$, $\mathcal{X} \in \mathbb{C}^{I_1 \times I_2 \times I_3}$ is a third-order tensor containing entries $x_{i_1, i_2, i_3}$'s.

Without loss of generality, one can partition the dimensions of a tensor into two groups, say $M$ and $N$ dimensions, separately. Thus, for two order-($M$+$N$) tensors: $\mathcal{X} \define (x_{i_1, \cdots, i_M, j_1, \cdots,j_N}) \in \mathbb{C}^{I_1 \times \cdots \times I_M\times
J_1 \times \cdots \times J_N}$ and $\mathcal{Y} \define (y_{i_1, \cdots, i_M, j_1, \cdots,j_N}) \in \mathbb{C}^{I_1 \times \cdots \times I_M\times
J_1 \times \cdots \times J_N}$, according to~\cite{MR3913666}, the \emph{tensor addition} $\mathcal{X} + \mathcal{Y}\in \mathbb{C}^{I_1 \times \cdots \times I_M\times
J_1 \times \cdots \times J_N}$ is given by 
\begin{eqnarray}\label{eq: tensor addition definition}
(\mathcal{X} + \mathcal{Y} )_{i_1, \cdots, i_M, j_1 \times \cdots \times j_N} &\define&
x_{i_1, \cdots, i_M, j_1 \times \cdots \times j_N} \nonumber \\
& &+ y_{i_1, \cdots, i_M, j_1 \times \cdots \times j_N}. 
\end{eqnarray}
On the other hand, for tensors $\mathcal{X} \define (x_{i_1, \cdots, i_M, j_1, \cdots,j_N}) \in \mathbb{C}^{I_1 \times \cdots \times I_M\times
J_1 \times \cdots \times J_N}$ and $\mathcal{Y} \define (y_{j_1, \cdots, j_N, k_1, \cdots,k_L}) \in \mathbb{C}^{J_1 \times \cdots \times J_N\times K_1 \times \cdots \times K_L}$, according to~\cite{MR3913666}, the \emph{Einstein product} (or simply referred to as \emph{tensor product} in this work) $\mathcal{X} \star_{N} \mathcal{Y} \in  \mathbb{C}^{I_1 \times \cdots \times I_M\times
K_1 \times \cdots \times K_L}$ is given by 
\begin{eqnarray}\label{eq: Einstein product definition}
\lefteqn{(\mathcal{X} \star_{N} \mathcal{Y} )_{i_1, \cdots, i_M,k_1 \times \cdots \times k_L} \define} \nonumber \\ &&\sum\limits_{j_1, \cdots, j_N} x_{i_1, \cdots, i_M, j_1, \cdots,j_N}y_{j_1, \cdots, j_N, k_1, \cdots,k_L}. 
\end{eqnarray}
Note that we will often abbreviate a tensor product $\mathcal{X} \star_{N} \mathcal{Y}$ to ``$\mathcal{X} \hspace{0.05cm}\mathcal{Y}$'' for notational simplicity in the rest of the paper. 
This tensor product will be reduced to the standard matrix multiplication as $L$ $=$ $M$ $=$ $N$ $=$ $1$. Other simplified situations can also be extended as tensor–vector product ($M >1$, $N=1$, and $L=0$) and tensor–matrix product ($M>1$ and $N=L=1$). In analogy to matrix analysis, we define some basic tensors and elementary tensor operations as follows. 

\begin{definition}\label{def: zero tensor}
A tensor whose entries are all zero is called a \emph{zero tensor}, denoted by $\mathcal{O}$. 
\end{definition}

\begin{definition}\label{def: identity tensor}
An \emph{identity tensor} $\mathcal{I} \in  \mathbb{C}^{I_1 \times \cdots \times I_N\times
J_1 \times \cdots \times J_N}$ is defined by 
\begin{eqnarray}\label{eq: identity tensor definition}
(\mathcal{I})_{i_1 \times \cdots \times i_N\times
j_1 \times \cdots \times j_N} \define \prod_{k = 1}^{N} \delta_{i_k, j_k},
\end{eqnarray}
where $\delta_{i_k, j_k} \define 1$ if $i_k  = j_k$; otherwise $\delta_{i_k, j_k} \define 0$.
\end{definition}
In order to define \emph{Hermitian} tensor, the \emph{conjugate transpose operation} (or \emph{Hermitian adjoint}) of a tensor is specified as follows.  
\begin{definition}\label{def: tensor conjugate transpose}
Given a tensor $\mathcal{X} \define (x_{i_1, \cdots, i_M, j_1, \cdots,j_N}) \in \mathbb{C}^{I_1 \times \cdots \times I_M\times J_1 \times \cdots \times J_N}$, its conjugate transpose, denoted by
$\mathcal{X}^{H}$, is defined by
\begin{eqnarray}\label{eq:tensor conjugate transpose definition}
(\mathcal{X}^H)_{ j_1, \cdots,j_N,i_1, \cdots, i_M}  \define  
\overline{x_{i_1, \cdots, i_M,j_1, \cdots,j_N}},
\end{eqnarray}
where the overline notion indicates the complex conjugate of the complex number $x_{i_1, \cdots, i_M,j_1, \cdots,j_N}$. If a tensor $\mathcal{X}$ satisfies $ \mathcal{X}^H = \mathcal{X}$, then $\mathcal{X}$ is a \emph{Hermitian tensor}. 
\end{definition}
\begin{definition}\label{def: unitary tensor}
Given a tensor $\mathcal{U} \define (u_{i_1, \cdots, i_M, i_1, \cdots,i_M}) \in \mathbb{C}^{I_1 \times \cdots \times I_M\times I_1 \times \cdots \times I_M}$, if
\begin{eqnarray}\label{eq:unitary tensor definition}
\mathcal{U}^H \star_M \mathcal{U} = \mathcal{U} \star_M \mathcal{U}^H = \mathcal{I} \in \mathbb{C}^{I_1 \times \cdots \times I_M\times I_1 \times \cdots \times I_M},
\end{eqnarray}
then $\mathcal{U}$ is a \emph{unitary tensor}. In this work, the symbol $\mathcal{U}$ is reserved for a unitary tensor. 
\end{definition}

Following definition is provided to define the inverse of a given tensor.
\begin{definition}\label{def: inverse of a tensor}
Given a \emph{square tensor} $\mathcal{X} \define (x_{i_1, \cdots, i_M, j_1, \cdots,j_M}) \in \mathbb{C}^{I_1 \times \cdots \times I_M\times I_1 \times \cdots \times I_M}$, if there exists $\mathcal{X} \in \mathbb{C}^{I_1 \times \cdots \times I_M\times I_1 \times \cdots \times I_M}$ such that 
\begin{eqnarray}\label{eq:tensor invertible definition}
\mathcal{X} \star_M \mathcal{X} = \mathcal{X} \star_M \mathcal{X} = \mathcal{I},
\end{eqnarray}
then $\mathcal{X}$ is the \emph{inverse} of $\mathcal{X}$. We usually write $\mathcal{X} \define \mathcal{X}^{-1}$ thereby. 
\end{definition}

We also list other crucial tensor operations here. The \emph{trace} of a square tensor is equivalent to the summation of all diagonal entries such that 
\begin{eqnarray}\label{eq: tensor trace def}
\mathrm{Tr}(\mathcal{X}) \define \sum\limits_{1 \leq i_j \leq I_j,\hspace{0.05cm}j \in [M]} \mathcal{X}_{i_1, \cdots, i_M,i_1, \cdots, i_M}.
\end{eqnarray}
The \emph{inner product} of two tensors $\mathcal{X}$, $\mathcal{Y} \in \mathbb{C}^{I_1 \times \cdots \times I_M\times J_1 \times \cdots \times J_N}$ is given by 
\begin{eqnarray}\label{eq: tensor inner product def}
\langle \mathcal{X}, \mathcal{Y} \rangle \define \mathrm{Tr}\left(\mathcal{X}^H \star_M \mathcal{Y}\right).
\end{eqnarray}
According to Eq.~\eqref{eq: tensor inner product def}, the \emph{Frobenius norm} of a tensor $\mathcal{X}$ is defined by 
\begin{eqnarray}\label{eq:Frobenius norm}
\left\Vert \mathcal{X} \right\Vert \define \sqrt{\langle \mathcal{X}, \mathcal{X} \rangle}.
\end{eqnarray}

%
%

From Theorem 5.2 in~\cite{ni2019hermitian}, every Hermitian tensor $\mathcal{H} \in  \mathbb{C}^{I_1 \times \cdots \times I_N \times I_1 \times \cdots \times I_N}$ has following decomposition
\begin{eqnarray}\label{eq:Hermitian Eigen Decom}
\mathcal{H} &=& \sum\limits_{i=1}^r \lambda_i \mathcal{U}_i \otimes \mathcal{U}^{H}_i, \mbox{
~with~~$\langle \mathcal{U}_i, \mathcal{U}_i \rangle =1$ and $\langle \mathcal{U}_i, \mathcal{U}_j \rangle = 0$ for $i \neq j$,}
\end{eqnarray}
where $\lambda_i \in \mathbb{R}$ and $\otimes$ denotes for Kronecker product. The values $\lambda_i$ are named as \emph{Hermitian eigevalues}, and the minimum integer of $r$ to decompose a Hermitian tensor as in Eq.~\eqref{eq:Hermitian Eigen Decom} is called \emph{Hermitian tensor rank}. A \emph{postive Hermitian tensor} is a Hermitian tensor with all \emph{Hermitian eigevalues} are positive. A \emph{nonnegative Hermitian tensor} is a Hermitian tensor with all \emph{Hermitian eigevalues} are nonnegative. The \emph{Hermitian determinant}, denoted as $\det\nolimits_H(\mathcal{A})$,  is defined as the product of $\lambda_i$ of the tensor $\mathcal{A}$. 

\subsection{Unitarily Invariant Tensor Norms}\label{sec:Unitarily Invariant Tensor Norms}


Let us represent the Hermitian eigenvalues of a Hermitian tensor $\mathcal{H} \in \mathbb{C}^{I_1 \times \cdots \times I_N \times I_1 \times \cdots \times I_N} $ in decreasing order by the vector $\vec{\lambda}(\mathcal{H}) = (\lambda_1(\mathcal{H}), \cdots, \lambda_r(\mathcal{H}))$, where $r$ is the Hermitian rank of the tensor $\mathcal{H}$. We use $\mathbb{R}_{\geq 0} (\mathbb{R}_{> 0})$ to represent a set of nonnegative (positive) real numbers. Let $\left\Vert \cdot \right\Vert_{\rho}$ be a unitarily invariant tensor norm, i.e., $\left\Vert \mathcal{H}\star_N \mathcal{U}\right\Vert_{\rho} = \left\Vert \mathcal{U}\star_N \mathcal{H}\right\Vert_{\rho} = \left\Vert \mathcal{H}\right\Vert_{\rho} $,  where $\mathcal{U}$ is any unitary tensor. Let $\rho : \mathbb{R}_{\geq 0}^r \rightarrow \mathbb{R}_{\geq 0}$ be the corresponding gauge function that satisfies H$\ddot{o}$lder’s inequality so that 
\begin{eqnarray}\label{eq:def gauge func and general unitarily inv norm}
\left\Vert \mathcal{H} \right\Vert_{\rho} = \left\Vert |\mathcal{H}| \right\Vert_{\rho} = \rho(\vec{\lambda}( | \mathcal{H} | ) ),
\end{eqnarray}
where $ |\mathcal{H}|  \define \sqrt{\mathcal{H}^H \star_N \mathcal{H}} $. The bijective correspondence between symmetric gauge functions on $\mathbb{R}_{\geq 0}^r$ and unitarily invariant norms is due to von Neumann~\cite{fan1955some}. 

Several popular norms can be treated as special cases of unitarily invariant tensor norm. The first one is Ky Fan like $k$-norm~\cite{fan1955some} for tensors. For $k \in \{1,2,\cdots, r \}$, the Ky Fan $k$-norm~\cite{fan1955some} for tensors  $\mathcal{H}  \mathbb{C}^{I_1 \times \cdots \times I_N \times I_1 \times \cdots \times I_N} $, denoted as $\left\Vert \mathcal{H}\right\Vert_{(k)}$, is defined as:
\begin{eqnarray}\label{eq:Ky Fan k norm for tensors}
\left\Vert \mathcal{H}\right\Vert_{(k)} \define \sum\limits_{i=1}^{k} \lambda_i(  |\mathcal{H}|  ).
\end{eqnarray}
If $k=1$,  the Ky Fan $k$-norm for tensors is the tensor operator norm (or spectral norm), denoted as $ \left\Vert \mathcal{H} \right\Vert$. The second one is Schatten $p$-norm for tensors, denoted as $\left\Vert \mathcal{H}\right\Vert_{p}$, is defined as:
\begin{eqnarray}\label{eq: Schatten p norm for tensors}
\left\Vert \mathcal{H}\right\Vert_{p} \define (\mathrm{Tr}|\mathcal{H}|^p )^{\frac{1}{p}},
\end{eqnarray}
where $ p \geq 1$. If $p=1$, it is the trace norm. The third one is $k$-trace norm, denoted as $\mathrm{Tr}_k[\mathcal{H}]$, defined by ~\cite{huang2020generalizing}. It is 
\begin{eqnarray}\label{eq: de k-trace norm for tensors}
\mathrm{Tr}_k[\mathcal{H}] \define \sum\limits_{1 \leq i_1 < i_2 < \cdots i_k \leq r} \lambda_{i_1} \lambda_{i_1}  \cdots \lambda_{i_k} 
\end{eqnarray}
where $ 1 \leq k \leq r$. If $k=1$, $\mathrm{Tr}_k[\mathcal{H}]$ is reduced as trace norm. In our later derivation for multivariate tensor norm inequalities in Section~\ref{sec:Multivariate Tensor Norm Inequalities} and tensor expander bounds in Section~\ref{sec:Tensor Expander Chernoff Bounds Derivation by Majorization}, we will focus on Ky Fan $k$-norm. 

Following inequality is the extension of H\"{o}lder inequality to gauge function $\rho$ which will be used by later majorization proof arguments. 
\begin{lemma}\label{lma:Holder inquality for gauge function}
For $n$ nonnegative real vectors with the dimension $r$, i.e., $\mathbf{b}_i = (b_{i_1}, \cdots, b_{i_r}) \in \mathbb{R}_{\geq 0}^r$, and $\alpha > 0$ with $\sum\limits_{i=1}^n \alpha_i = 1$, we have 
\begin{eqnarray}\label{eq1:lma:Holder inquality for gauge function}
\rho\left( \prod\limits_{i=1}^n b_{i_1}^{\alpha_i},  \prod\limits_{i=1}^n b_{i_2}^{\alpha_i}, \cdots,  \prod\limits_{i=1}^n b_{i_r}^{\alpha_i}  \right) \leq  \prod\limits_{i=1}^n \rho^{\alpha_i} (\mathbf{b}_i)
\end{eqnarray}
\end{lemma}
\textbf{Proof:}
This proof is based on mathematical induction. The base case for $n=2$ has been shown by Theorem IV.1.6 from~\cite{bhatia2013matrix}. 

We assume that Eq.~\eqref{eq1:lma:Holder inquality for gauge function} is true for $n=m$, where $m > 2$. Let $\odot$ be the component-wise product (Hadamard product) between two vectors.  Then, we have 
\begin{eqnarray}\label{eq2:lma:Holder inquality for gauge function}
\rho\left( \prod\limits_{i=1}^{m+1} b_{i_1}^{\alpha_i},  \prod\limits_{i=1}^{m+1} b_{i_2}^{\alpha_i}, \cdots,  \prod\limits_{i=1}^{m+1} b_{i_r}^{\alpha_i}  \right) = 
\rho\left( \odot_{i=1}^{m+1} \mathbf{b}_i^{\alpha_i}  \right),
\end{eqnarray}
where $\odot_{i=1}^{m+1} \mathbf{b}_i^{\alpha_i}$ is defined as $\left( \prod\limits_{i=1}^{m+1} b_{i_1}^{\alpha_i},  \prod\limits_{i=1}^{m+1} b_{i_2}^{\alpha_i}, \cdots,  \prod\limits_{i=1}^{m+1} b_{i_r}^{\alpha_i}  \right)$ with $\mathbf{b}_i^{\alpha_i} \define (b_{i_1}^{\alpha_i}, \cdots, b_{i_r}^{\alpha_i})$. Under such notations, Eq.~\eqref{eq2:lma:Holder inquality for gauge function} can be bounded as  
\begin{eqnarray}\label{eq3:lma:Holder inquality for gauge function}
\rho\left( \odot_{i=1}^{m+1} \mathbf{b}_i^{\alpha_i}  \right) &= &
\rho\left( \left( \odot_{i=1}^{m} \mathbf{b}_i^{\frac{\alpha_i}{ \sum\limits_{j=1}^m \alpha_j }  } \right)^{\sum\limits_{j=1}^m \alpha_j } \odot \mathbf{b}_{m+1}^{\alpha_{m+1}}\right) \nonumber \\
& \leq & \left[ \rho^{\sum\limits_{j=1}^m \alpha_j } \left( \odot_{i=1}^{m} \mathbf{b}_i^{\frac{\alpha_i}{ \sum\limits_{j=1}^m \alpha_j }  }  \right)  \right] \cdot \rho^{\alpha_{m+1}}( \mathbf{b}_{m+1}) \leq  \prod\limits_{i=1}^{m+1} \rho^{\alpha_i}(\mathbf{b}_i). 
\end{eqnarray}
By mathematical induction, this lemma is proved. $\hfill \Box$

\subsection{Antisymmetric Tensor Product}\label{sec:Antisymmetric Tensor Product}


Let $\mathfrak{H}$ be a space of Hermitian tensors with Hermitian rank $r$. Two tensors $\mathcal{C}, \mathcal{B} \in \mathfrak{H}$ is said $\mathcal{C} \geq \mathcal{B}$ if $\mathcal{C} - \mathcal{B} $ is a nonnegative Hermitian tensor. For any $k \in \{1,2,\cdots,r\}$, let $\mathfrak{H}^{\otimes k}$ be the $k$-th tensor power of the tensor space $\mathfrak{H}$ and let $\mathfrak{H}^{\wedge k}$ be the antisymmetric subspace of $\mathfrak{H}^{\otimes k}$. The $k$-th antisymmetric tensor power, $\wedge^k : \mathcal{C} \rightarrow \mathcal{C}^{\wedge k}$, maps any Hermitian tensor $\mathcal{C}$ to the restriction of $\mathcal{C}^{\otimes k} \in \mathfrak{H}^{\otimes k}$ to the antisymmetric subspace $\mathfrak{H}^{\wedge k}$ of $\mathfrak{H}^{\otimes k}$. Following lemma summarizes several useful properties of such antisymmetric tensor products. 

\begin{lemma}\label{lma:antisymmetric tensor product properties}
Let $\mathcal{C}, \mathcal{B}, \mathcal{C}, \mathcal{D} \in \mathbb{C}^{I_1 \times \cdots \times I_N \times I_1 \times \cdots \times I_N}$ be Hermitian tensors from $\mathfrak{H}$ with Hermitian rank $r$. For any $k \in \{1,2,\cdots,r\}$, we have 
\begin{enumerate}[label={[\arabic*]}]
	\item $(\mathcal{C}^{\wedge k})^H = (\mathcal{C}^H)^{\wedge k}$,.
	\item $(\mathcal{C}^{\wedge k}) \star_N (\mathcal{B}^{\wedge k})= (\mathcal{C}\star_N \mathcal{B})^{\wedge k}$. 
	\item If $\lim\limits_{i \rightarrow \infty} \left\Vert \mathcal{C}_i -  \mathcal{C} \right\Vert \rightarrow 0$ , then $\lim\limits_{i \rightarrow \infty} \left\Vert \mathcal{C}^{\wedge k}_i -  \mathcal{C}^{\wedge k} \right\Vert \rightarrow 0$.
	\item If $\mathcal{C} \geq \mathcal{O}$ (zero tensor), then $\mathcal{C}^{\wedge k} \geq \mathcal{O}$ and $(\mathcal{C}^p)^{\wedge k} = (\mathcal{C}^{\wedge k})^p$ for all $p \in \mathbb{R}_{\ge 0 }$.
    \item  $|\mathcal{C}|^{\wedge k} = | \mathcal{C}^{\wedge k}|$.
    \item  If $\mathcal{D} \geq \mathcal{O}$ and $\mathcal{D}$ is invertibale,  $(\mathcal{D}^z)^{\wedge k} = (\mathcal{D}^{\wedge k})^z$ for all $z \in \mathbb{C}$.
    \item  $\left\Vert \mathcal{C}^{\wedge k} \right\Vert = \prod\limits_{i=1}^{k} \lambda_i ( | \mathcal{C} |)$.
\end{enumerate}
\end{lemma}
\textbf{Proof:}
Facts $\textit{[1]}$ and $\textit{[2]}$ are the restrictions of the associated relations $(\mathcal{C}^H)^{\otimes k} = (\mathcal{C}^{\otimes k})^H$ and $(\mathcal{C} \star_N \mathcal{B})^{\otimes k} = (\mathcal{C}^{\otimes k})\star_N (\mathcal{B}^{\otimes k})$ to $\mathfrak{H}^{\wedge k}$. The fact $\textit{[3]}$ is true since, if $\lim\limits_{i \rightarrow \infty} \left\Vert \mathcal{C}_i -  \mathcal{C} \right\Vert \rightarrow 0$, we have $\lim\limits_{i \rightarrow \infty} \left\Vert \mathcal{C}^{\otimes k}_i -  \mathcal{C}^{\otimes k} \right\Vert \rightarrow 0$ and the asscoaited restrictions of $\mathcal{C}_i^{\otimes k}, \mathcal{C}^{\otimes k}$ to the antisymmetric subspace $\mathfrak{H}^k$. 

For the fact $\textit{[4]}$, if $\mathcal{C} \geq \mathcal{O}$, then we have $\mathcal{C}^{\wedge k} = ((\mathcal{C}^{1/2})^{\wedge k})^H \star_N  ((\mathcal{C}^{1/2})^{\wedge k}) \geq   \mathcal{O}$ from facts  $\textit{[1]}$ and $\textit{[2]}$. If $p$ is ratonal, we have  $(\mathcal{C}^p)^{\wedge k} = (\mathcal{C}^{\wedge k})^p$  from the fact $\textit{[2]}$, and the equality $(\mathcal{C}^p)^{\wedge k} = (\mathcal{C}^{\wedge k})^p$ is also true for any $p > 0$ if we apply the fact $\textit{[3]}$ to approximate any irrelational numbers by rational numbers.

Because we have 
\begin{eqnarray}
|\mathcal{C}|^{\wedge k} =  \left( \sqrt{\mathcal{C}^H \mathcal{C}}\right)^{\wedge k}  =   \sqrt{ (\mathcal{C}^{\wedge k})^H \mathcal{C}^{\wedge k}  }=  | \mathcal{C}^{\wedge k}|,
\end{eqnarray}
from facts $\textit{[1]}$, $\textit{[2]}$ and $\textit{[4]}$, so the fact $\textit{[5]}$ is valid. 

For the fact $\textit{[6]}$, if $z  < 0$, the fact $\textit{[6]}$ is true for all $z \in \mathbb{R}$ by applying the fact $\textit{[4]}$ to $\mathcal{D}^{-1}$. Since we can apply the definition $\mathcal{D}^{z} \define \exp(z \ln \mathcal{D})$ to have
\begin{eqnarray}
\mathcal{C}^p &=& \mathcal{D}^z~~\leftrightarrow~~\mathcal{C} = \exp\left(\frac{z}{p} \ln \mathcal{D} \right),
\end{eqnarray}
where $\mathcal{C} \geq \mathcal{O}$. The general case of any $z \in \mathbb{C}$ is also true by applying the fact $\textit{[4]}$ to $\mathcal{C} = \exp(\frac{z}{p} \ln \mathcal{D})$. 

For the fact $\textit{[7]}$ proof, it is enough to prove the case that $\mathcal{C} \geq \mathcal{O}$ due to the fact $\textit{[5]}$. Then, there exisits a set of orthogonal tensors $\{\mathcal{U}_1, \cdots, \mathcal{U}_r\}$ such that $\mathcal{C} \star_N \mathcal{U}_i = \lambda_i   \mathcal{U}_i$ for $1 \leq i \leq r$. We then have 
\begin{eqnarray}
\mathcal{C}^{\wedge k} \left(\mathcal{U}_{i_1} \wedge \cdots \wedge \mathcal{U}_{i_k}\right)
&=& \mathcal{C}\star_N  \mathcal{U}_{i_1} \wedge \cdots \wedge \mathcal{C}\star_N  \mathcal{U}_{i_k}  \nonumber \\
&=& \left( \prod\limits_{i=1}^{k} \lambda_i ( | \mathcal{C} |)  \right) \mathcal{U}_{i_1} \wedge \cdots \wedge \mathcal{U}_{i_k}.
\end{eqnarray}
Hence, $\left\Vert \mathcal{C}^{\wedge k} \right\Vert = \prod\limits_{i=1}^{k} \lambda_i ( | \mathcal{C} |)$. $\hfill \Box$

\subsection{Majorization}\label{sec:Majorization}


In this subsection, we will discuss majorization and several lemmas about majorization which will be used at later proofs. 

Let $\mathbf{x} = [x_1, \cdots,x_r] \in \mathbb{R}^{r}, \mathbf{y} = [y_1, \cdots,y_r] \in \mathbb{R}^{r}$ be two vectors with following orders among entries $x_1 \geq \cdots \geq x_r$ and $y_1 \geq \cdots \geq y_r$, \emph{weak majorization} between vectors $\mathbf{x}, \mathbf{y}$, represented by $\mathbf{x} \prec_{w} \mathbf{y}$, requires following relation for  vectors $\mathbf{x}, \mathbf{y}$:
\begin{eqnarray}\label{eq:weak majorization def}
\sum\limits_{i=1}^k x_i \leq \sum\limits_{i=1}^k y_i,
\end{eqnarray}
where $k \in \{1,2,\cdots,r\}$. \emph{Majorization} between vectors $\mathbf{x}, \mathbf{y}$, indicated by $\mathbf{x} \prec \mathbf{y}$, requires following relation for vectors $\mathbf{x}, \mathbf{y}$:
\begin{eqnarray}\label{eq:majorization def}
\sum\limits_{i=1}^k x_i &\leq& \sum\limits_{i=1}^k y_i,~~\mbox{for $1 \leq k < r$;} \nonumber \\
\sum\limits_{i=1}^r x_i &=& \sum\limits_{i=1}^r y_i,~~\mbox{for $k = r$.}
\end{eqnarray}

For $\mathbf{x}, \mathbf{y} \in \mathbb{R}^r_{\geq 0}$ such that  $x_1 \geq \cdots \geq x_r$ and $y_1 \geq \cdots \geq y_r$,  \emph{weak log majorization} between vectors $\mathbf{x}, \mathbf{y}$, represented by $\mathbf{x} \prec_{w \log} \mathbf{y}$, requires following relation for vectors $\mathbf{x}, \mathbf{y}$:
\begin{eqnarray}\label{eq:weak log majorization def}
\prod\limits_{i=1}^k x_i \leq \prod\limits_{i=1}^k y_i,
\end{eqnarray}
where $k \in \{1,2,\cdots,r\}$, and \emph{log majorization} between vectors $\mathbf{x}, \mathbf{y}$, represented by $\mathbf{x} \prec_{\log} \mathbf{y}$, requires equality for $k=r$ in Eq.~\eqref{eq:weak log majorization def}. If $f$ is a single variable function, $f(\mathbf{x})$ represents a vector of $[f(x_1),\cdots,f(x_r)]$. From Lemma 1 in~\cite{hiai2017generalized}, we have 
\begin{lemma}\label{lma:Lemma 1 Gen Log Hiai}
(1) For any convex function $f: [0, \infty) \rightarrow [0, \infty)$, if we have $\mathbf{x} \prec \mathbf{y}$, then $f(\mathbf{x}) \prec_{w} f(\mathbf{y})$. \\
(2) For any convex function and non-decreasing $f: [0, \infty) \rightarrow [0, \infty)$, if we have $\mathbf{x} \prec_{w} \mathbf{y}$, then $f(\mathbf{x}) \prec_{w} f(\mathbf{y})$. \\
\end{lemma}

Another lemma is from Lemma 12 in~\cite{hiai2017generalized}, we have 
\begin{lemma}\label{lma:Lemma 12 Gen Log Hiai}
Let $\mathbf{x}, \mathbf{y} \in \mathbb{R}^r_{\geq 0}$ such that  $x_1 \geq \cdots \geq x_r$ and $y_1 \geq \cdots \geq y_r$ with $\mathbf{x}\prec_{\log} \mathbf{y}$. Also let $\mathbf{y}_i = [y_{i;1}, \cdots , y_{i;r} ] \in \mathbb{R}^r_{\geq 0}$ be a sequence of vectors such that $y_{i;1} \geq \cdots \geq y_{i;r} > 0$ and $\mathbf{y}_i \rightarrow \mathbf{y}$ as $i \rightarrow \infty$. Then, there exists $i_0 \in \mathbb{N}$ and $\mathbf{x}_i  = [x_{i;1}, \cdots , x_{i;r} ] \in \mathbb{R}^r_{\geq 0}$ for $i \geq i_0$ such that $x_{i;1} \geq \cdots \geq x_{i;r} > 0$, $\mathbf{x}_i \rightarrow \mathbf{x}$ as $i \rightarrow \infty$, and 
\begin{eqnarray}
\mathbf{x}_i \prec_{\log} \mathbf{y}_i \mbox{~~for $i \geq i_0$.} 
\end{eqnarray}
\end{lemma}

For any function $f$ on $\mathbb{R}_{\geq 0}$, the term $f(\mathbf{x})$ is defined as $f(\mathbf{x}) \define (f(x_1), \cdots, f(x_r))$ with conventions $e^{ - \infty} = 0$ and $\log 0 = - \infty$. 

\section{Multivariate Tensor Norm Inequalities}\label{sec:Multivariate Tensor Norm Inequalities}

In this section, we will develop several theorems about majorization in Section~\ref{sec:Majorization wtih Integral Average}, and log majorization with integral average in Section~\ref{sec:Log-Majorization wtih Integral Average}. These majorization related theorems will provide us tools  in deriving bounds for Ky Fan $k$-norms of multivariate tensors in Section~\ref{sec:Multivaraite Tensor Norm Inequalities}.

\subsection{Majorization wtih Integral Average}\label{sec:Majorization wtih Integral Average}


Let $\Omega$ be a $\sigma$-compact metric space and $\nu$ a probability measure on the Borel $\sigma$-field of $\Omega$. Let $\mathcal{C}, \mathcal{D}_\tau \in \mathbb{C}^{I_1 \times \cdots \times I_N \times I_1 \times \cdots \times I_N}$ be Hermitian tensors with Hermitian rank $r$. We further assume that tensors $\mathcal{C}, \mathcal{D}_\tau$ are uniformly bounded in their norm for $\tau \in \Omega$. Let $\tau \in\Omega \rightarrow  \mathcal{D}_\tau$ be a continuous function such that $\sup \{\left\Vert  D_{\tau} \right\Vert: \tau \in \Omega  \} < \infty$. For notational convenience, we define the following relation:
\begin{eqnarray}\label{eq:integral eigen vector rep}
\left[ \int_{\Omega} \lambda_1(\mathcal{D}_\tau) d\nu(\tau), \cdots, \int_{\Omega} \lambda_r(\mathcal{D}_\tau) d\nu(\tau) \right] \define \int_{\Omega^r} \vec{\lambda}(\mathcal{D}_\tau) d\nu^r(\tau).
\end{eqnarray}
If $f$ is a single variable function, the notation $f(\mathcal{C})$ represents a tensor function with respect to the tensor $\mathcal{C}$. We want to prove the following theorem about weak majorization of the integral average. 

\begin{theorem}\label{thm:weak int average thm 4}
Let $\Omega, \nu, \mathcal{C}, \mathcal{D}_\tau$ be defined as the beginning part of Section~\ref{sec:Majorization wtih Integral Average}, and $f: \mathbb{R} \rightarrow [0, \infty)$ be a non-decreasing convex function, we have following two equivalent statements:
\begin{eqnarray}\label{eq1:thm:weak int average thm 4}
\vec{\lambda}(\mathcal{C}) \prec_w  \int_{\Omega^r} \vec{\lambda}(\mathcal{D}_\tau) d\nu^r(\tau) \Longleftrightarrow \left\Vert f(\mathcal{C}) \right\Vert_{(k)} \leq 
\int_{\Omega} \left\Vert f(\mathcal{D}_{\tau}) \right\Vert_{(k)}  d\nu(\tau),
\end{eqnarray}
where $\left\Vert \cdot \right\Vert_{(k)}$ is the Ky Fan $k$-norm defined by Eq.~\eqref{eq:Ky Fan k norm for tensors}.
\end{theorem}
\textbf{Proof:}
We assume that the left statement of Eq.~\eqref{eq1:thm:weak int average thm 4} is true and the function $f$ is a non-decreasing convex function. From Lemma~\ref{lma:Lemma 1 Gen Log Hiai}, we have 
\begin{eqnarray}\label{eq2:thm:weak int average thm 4}
\vec{\lambda}(f (\mathcal{C})) = f (\vec{\lambda}(\mathcal{C})) \prec_w  f \left(\int_{\Omega^r} \vec{\lambda}(\mathcal{D}_\tau) d\nu^r(\tau) \right).
\end{eqnarray}
From the convexity of $f$, we also have 
\begin{eqnarray}\label{eq3:thm:weak int average thm 4}
f \left(\int_{\Omega^r} \vec{\lambda}(\mathcal{D}_\tau) d\nu^r(\tau) \right) \leq \int_{\Omega^r} f(\vec{\lambda}(\mathcal{D}_\tau)) d\nu^r(\tau) = \int_{\Omega^r} \vec{\lambda} ( f(\mathcal{D}_\tau)) d\nu^r(\tau).
\end{eqnarray}
Then, we obtain $\vec{\lambda}(f (\mathcal{C}))  \prec_{w} \int_{\Omega^r} \vec{\lambda} ( f(\mathcal{D}_\tau)) d\nu^r(\tau)$. By applying Lemma 4.4.2 in~\cite{hiai2010matrix} to both sides of $\vec{\lambda}(f (\mathcal{C}))  \prec_{w}  \int_{\Omega^r} \vec{\lambda} ( f(\mathcal{D}_\tau)) d\nu^r(\tau)$ with gauge function $\rho$ of Ky Fan $k$-norm, we obtain 
\begin{eqnarray}\label{eq4:thm:weak int average thm 4}
\left \Vert f(\mathcal{C}) \right\Vert_{(k)} &\leq &\rho \left( \int_{\Omega^r} \vec{\lambda} ( f(\mathcal{D}_\tau)) d\nu^r(\tau)  \right)  \nonumber \\
&\leq & \int_{\Omega} \rho(\vec{\lambda} ( f(\mathcal{D}_\tau))) d\nu(\tau) 
= \int_{\Omega} \left\Vert f(\mathcal{D}_\tau) \right\Vert_{(k)} d\nu(\tau).
\end{eqnarray}
Therefore, the right statement of Eq.~\eqref{eq1:thm:weak int average thm 4} is true from the left statement. 

On the other hand, if the right statement of Eq.~\eqref{eq1:thm:weak int average thm 4} is true, we select a function $f \define \max\{x + c, 0\} $, where $c$ is a postive real constant satisfying $\mathcal{C} + c \mathcal{I} \geq \mathcal{O}$, $\mathcal{D}_{\tau} + c \mathcal{I} \geq \mathcal{O}$ for all $\tau \in \Omega$, and tensors $\mathcal{C} + c \mathcal{I}, \mathcal{D}_{\tau} + c \mathcal{I}$ having Hermitian rank $r$. If the Ky Fan norm $k$-norm at the right statement of Eq.~\eqref{eq1:thm:weak int average thm 4} is applied, we have 
\begin{eqnarray}\label{eq5:thm:weak int average thm 4}
\sum\limits_{i=1}^k (\lambda_i (\mathcal{C}) + c ) \leq  
\sum\limits_{i=1}^k \int_{\Omega} ( \lambda_i (\mathcal{D}_{\tau}) + c ) d\nu(\tau).
\end{eqnarray}
Hence, $\sum\limits_{i=1}^k \lambda_i (\mathcal{C}) \leq  
\sum\limits_{i=1}^k \int_{\Omega} \lambda_i (\mathcal{D}_{\tau}) d\nu(\tau)$, this is the left statement of Eq.~\eqref{eq1:thm:weak int average thm 4}. $\hfill \Box$

Next theorem will provide a stronger version of Theorem~\ref{thm:weak int average thm 4} by enhancing weak majorization to majorization. 
\begin{theorem}\label{thm:weak int average thm 5}
Let $\Omega, \nu, \mathcal{C}, \mathcal{D}_\tau$ be defined as the beginning part of Section~\ref{sec:Majorization wtih Integral Average}, and $f: \mathbb{R} \rightarrow [0, \infty)$ be a convex function, we have following two equivalent statements:
\begin{eqnarray}\label{eq1:thm:weak int average thm 5}
\vec{\lambda}(\mathcal{C}) \prec  \int_{\Omega^r} \vec{\lambda}(\mathcal{D}_\tau) d\nu^r(\tau) \Longleftrightarrow \left\Vert f(\mathcal{C}) \right\Vert_{(k)} \leq 
\int_{\Omega} \left\Vert f(\mathcal{D}_{\tau}) \right\Vert_{(k)}  d\nu(\tau).
\end{eqnarray} 
\end{theorem}
\textbf{Proof:}
We assume that the left statement of Eq.~\eqref{eq1:thm:weak int average thm 5} is true and the function $f$ is a convex function. Again, from Lemma~\ref{lma:Lemma 1 Gen Log Hiai}, we have
\begin{eqnarray}\label{eq2:thm:weak int average thm 5}
\vec{\lambda}(f(\mathcal{C})) = f(\vec{\lambda}(\mathcal{C})) \prec_{w}  f \left( \left(\int_{\Omega^r} \vec{\lambda}(\mathcal{D}_\tau) d\nu^r(\tau)  \right) \right) \leq \int_{\Omega^r} f(\vec{\lambda}(\mathcal{D}_\tau)) d\nu^r(\tau),
\end{eqnarray}
then, 
\begin{eqnarray}\label{eq3:thm:weak int average thm 5}
\left\Vert f(\mathcal{C}) \right\Vert_{(k)} &\leq & \rho\left( \int_{\Omega^r} f(\vec{\lambda}(\mathcal{D}_\tau)) d\nu^r(\tau) \right) \nonumber \\
& \leq & \int_{\Omega}\rho \left( f(\vec{\lambda}(\mathcal{D}_\tau)) \right)d\nu (\tau) = 
\int_{\Omega} \left\Vert f( \mathcal{D}_\tau) \right\Vert_{(k)} d\nu (\tau),
\end{eqnarray}
where $\rho$ is the gauge function of Ky Fan $k$-norm. This proves the right statement of Eq.~\eqref{eq1:thm:weak int average thm 5}. 

Now, we assume that the right statement of Eq.~\eqref{eq1:thm:weak int average thm 5} is true. From Theorem~\ref{thm:weak int average thm 4}, we already have $\vec{\lambda}(\mathcal{C}) \prec_w  \int_{\Omega^r} \vec{\lambda}(\mathcal{D}_\tau) d\nu^r(\tau)$.  It is enough to prove $\sum\limits_{i=1}^r \lambda_i(\mathcal{C}) \geq \int_{\Omega} \sum\limits_{i=1}^r \lambda_i(\mathcal{D}_{\tau}) d \nu(\tau)$. We define a function $f \define \max\{c - x, 0\} $, where $c$ is a postive real constant satisfying $\mathcal{C} \leq c \mathcal{I} $, $\mathcal{D}_{\tau} \leq  c \mathcal{I}$ for all $\tau \in \Omega$ and tensors $c \mathcal{I} - \mathcal{C}, c \mathcal{I} - \mathcal{D}_{\tau}$ having Hermitian rank $r$. If the trace norm is applied, i.e., the sum of the absolute value of all eigenvalues of a Hermitian tensor, then the right statement of Eq.~\eqref{eq1:thm:weak int average thm 5} becomes
\begin{eqnarray}\label{eq4:thm:weak int average thm 5}
\sum\limits_{i=1}^r \lambda_i \left( c\mathcal{I} - \mathcal{C}\right)  \leq \int_{\Omega} 
\sum\limits_{i=1}^r \lambda_i \left( c\mathcal{I} - \mathcal{D}_{\tau}\right) d \nu(\tau).
\end{eqnarray}
The desired inequality  $\sum\limits_{i=1}^r \lambda_i(\mathcal{C}) \geq \int_{\Omega} \sum\limits_{i=1}^r \lambda_i(\mathcal{D}_{\tau}) d \nu(\tau)$ is established. $\hfill \Box$

\subsection{Log-Majorization wtih Integral Average}\label{sec:Log-Majorization wtih Integral Average}


The purpose of this section is to consider log-majorization issues for Ky Fan $k$-norm of Hermitian tensors. In this section, let $\mathcal{C}, \mathcal{D}_\tau \in \mathbb{C}^{I_1 \times \cdots \times I_N \times I_1 \times \cdots \times I_N}$ be nonnegative Hermitian tensors with Hermitian rank $r$, i,e, all Hermitian eigenvalues are positive, and keep other notations with the same definitions as at the beginning of the Section~\ref{sec:Majorization wtih Integral Average}. For notational convenience, we define the following relation for logarithm vector:
\begin{eqnarray}\label{eq:integral eigen log vector rep}
\left[ \int_{\Omega} \log \lambda_1(\mathcal{D}_\tau) d\nu(\tau), \cdots, \int_{\Omega} \log \lambda_r(\mathcal{D}_\tau) d\nu(\tau) \right] \define \int_{\Omega^r} \log \vec{\lambda}(\mathcal{D}_\tau) d\nu^r(\tau).
\end{eqnarray}

Following theorem is used to build the relationship between weak log-majorization of eigenvalues and Ky Fan $k$-norm.
\begin{theorem}\label{thm:weak int log average thm 7}
Let $\mathcal{C}, \mathcal{D}_\tau$ be nonnegative Hermitian tensors, $f: (0, \infty) \rightarrow [0,\infty)$ be a continous function such that the mapping $x \rightarrow \log f(e^x)$ is convex on $\mathbb{R}$, and $g: (0, \infty) \rightarrow [0,\infty)$ be a continous function such that the mapping $x \rightarrow g(e^x)$ is convex on $\mathbb{R}$ , then  we have following three equivalent statements:
\begin{eqnarray}\label{eq1:thm:weak int log average thm 7}
\vec{\lambda}(\mathcal{C}) &\prec_{w \log}& \exp  \int_{\Omega^r} \log \vec{\lambda}(\mathcal{D}_\tau) d\nu^r(\tau);
\end{eqnarray}
\begin{eqnarray}\label{eq2:thm:weak int log average thm 7}
\left\Vert f(\mathcal{C}) \right\Vert_{(k)} &\leq &
\exp \int_{\Omega} \log \left\Vert f(\mathcal{D}_{\tau}) \right\Vert_{(k)}  d\nu(\tau);
\end{eqnarray}
\begin{eqnarray}\label{eq3:thm:weak int log average thm 7}
\left\Vert g(\mathcal{C}) \right\Vert_{(k)} &\leq &
\int_{\Omega} \left\Vert g(\mathcal{D}_{\tau}) \right\Vert_{(k)}  d\nu(\tau).
\end{eqnarray}
\end{theorem}
\textbf{Proof:}
The roadmap of this proof is to prove equivalent statements between Eq.~\eqref{eq1:thm:weak int log average thm 7} and Eq.~\eqref{eq2:thm:weak int log average thm 7} first, followed by equivalent statements between Eq.~\eqref{eq1:thm:weak int log average thm 7} and Eq.~\eqref{eq3:thm:weak int log average thm 7}. 

\textbf{Eq.~\eqref{eq1:thm:weak int log average thm 7} $\Longrightarrow$ Eq.~\eqref{eq2:thm:weak int log average thm 7}}

There are two cases to be discussed in this part of proof: $\mathcal{C}, \mathcal{D}_\tau$ are positive Hermitian tensors, and $\mathcal{C}, \mathcal{D}_\tau$ are nonnegative Hermitian tensors. We consider the case that $\mathcal{C}, \mathcal{D}_\tau$ are positive Hermitian tensors first.

Since $\mathcal{D}_\tau$ are positive, we can find $\varepsilon > 0$ such that $\mathcal{D}_{\tau} \geq \varepsilon \mathcal{I}$ for all $\tau \in \Omega$. From Eq.~\eqref{eq1:thm:weak int log average thm 7}, the convexity of $\log f (e^x)$ and Lemma~\ref{lma:Lemma 1 Gen Log Hiai}, we have 
\begin{eqnarray}
\vec{\lambda} \left( f ( \mathcal{C})\right)  = f \left(\exp \left( \log \vec{\lambda} (\mathcal{C}) \right) \right) &\prec _w &  f \left(\exp    \int_{\Omega^r} \log \vec{\lambda}(\mathcal{D}_\tau) d\nu^r(\tau)             \right) \nonumber \\
& \leq &  \exp \left( \int_{\Omega^r} \log f \left( \vec{\lambda}(\mathcal{D}_\tau) \right)  d\nu^r(\tau) \right).
\end{eqnarray}
Then, from Eq.~\eqref{eq:def gauge func and general unitarily inv norm}, we obtain
\begin{eqnarray}\label{eq4:thm:weak int log average thm 7}
\left\Vert f (\mathcal{C}) \right\Vert _{\rho}
& \leq &  \rho\left(\exp \left( \int_{\Omega^r} \log  f \left( \vec{\lambda}(\mathcal{D}_\tau) \right) d\nu^r(\tau) \right)   \right).
\end{eqnarray}


From the function $f$ properties, we can assume that $f(x) > 0$ for any $x > 0$. Then, we have 
following bounded and continous maps on $\Omega$: $\tau \rightarrow \log f (\lambda_i (\mathcal{D}_{\tau}))  $ for $i \in \{1,2,\cdots, r\}$, and $\tau \rightarrow \log \left\Vert f (\mathcal{D}_{\tau}) \right\Vert_{(k)}$. Because we have $\nu (\Omega) = 1$ and $\sigma$-compactness of $\Omega$, we have $\tau_{k}^{(n)} \in \Omega$ and $\alpha_{k}^{(n)}$ for $k \in \{1,2,\cdots, n\}$ and $n \in \mathbb{N}$ with $\sum\limits_{k=1}^{n} \alpha_{k}^{(n)} = 1$ such that 
\begin{eqnarray}\label{eq:35}
\int_{\Omega} \log f (\lambda_i ( \mathcal{D}_{\tau} ) ) d \nu (\tau) = \lim\limits_{n \rightarrow \infty} \sum\limits_{k=1}^{n} \alpha_{k}^{(n)}  \log f (\lambda_i (\mathcal{D}_{\tau_k^{(n)} }))   , \mbox{for $i \in \{1,2,\cdots r \}$};
\end{eqnarray}
and 
\begin{eqnarray}\label{eq:36}
\int_{\Omega} \log \left\Vert f (\mathcal{D}_{\tau}) \right\Vert_{(k)} d \nu (\tau) = \lim\limits_{n \rightarrow \infty} \sum\limits_{k=1}^{n} \alpha_{k}^{(n)}  \log \left\Vert f (\mathcal{D}_{\tau_k^{(n)}    }) \right\Vert_{(k)} .
\end{eqnarray}
By taking the exponential at both sides of Eq.~\eqref{eq:35} and apply the gauge function $\rho$, we have
\begin{eqnarray}\label{eq:37}
\rho \left( \exp \int_{\Omega^r} \log f (\vec{\lambda} ( \mathcal{D}_{\tau} ) ) d \nu^r (\tau)  \right)= \lim\limits_{n \rightarrow \infty} \rho\left(  \prod\limits_{k=1}^{n}  f  \left( \vec{\lambda} \left(\mathcal{D}_{\tau_k^{(n)} } \right)  \right)^{ \alpha_{k}^{(n)} }  \right).
\end{eqnarray}
Similarly, by taking the exponential at both sides of Eq.~\eqref{eq:36}, we have
\begin{eqnarray}\label{eq:38}
\exp \left( \int_{\Omega} \log \left\Vert f (\mathcal{D}_{\tau}) \right\Vert_{(k)} d \nu (\tau) \right) = \lim\limits_{n \rightarrow \infty} \prod \limits_{k=1}^{n} \left\Vert f \left( \mathcal{D}_{\tau_k^{(n)}    } \right) \right\Vert^{\alpha_{k}^{(n)}}_{(k)} .
\end{eqnarray}
From Lemma~\ref{lma:Holder inquality for gauge function}, we have 
\begin{eqnarray}\label{eq:41}
\rho \left(  \prod\limits_{k=1}^{n}  f  \left( \vec{\lambda} \left(\mathcal{D}_{\tau_k^{(n)} } \right)  \right)^{ \alpha_{k}^{(n)} }  \right) & \leq & \prod \limits_{k=1}^{n} \rho^{  \alpha_{k}^{(n)}    } \left( f \left( \vec{\lambda} \left( \mathcal{D}_{ \tau_{k}^{(n)}  }\right) 
 \right) \right) \nonumber \\
&=& \prod \limits_{k=1}^{n} \rho^{  \alpha_{k}^{(n)}    } \left( \vec{\lambda}  \left( f  \left( \mathcal{D}_{ \tau_{k}^{(n)}  }\right) 
 \right) \right) \nonumber \\
&=& \prod \limits_{k=1}^{n} \left\Vert f  \left( \mathcal{D}_{ \tau_{k}^{(n)}  } \right) \right\Vert_{(k)}^{\alpha_{k}^{(n)}}
\end{eqnarray}

From Eqs.~\eqref{eq:37},~\eqref{eq:38} and~\eqref{eq:41}, we have 
\begin{eqnarray}\label{eq:42}
\rho \left( \exp \int_{\Omega^r} \log f (\vec{\lambda} ( \mathcal{D}_{\tau} ) ) d \nu^r (\tau)  \right) \leq \exp \int_{\Omega} \log \left\Vert f(\mathcal{D}_{\tau})\right\Vert_{(k)} d \nu(\tau).
\end{eqnarray}
Then, Eq.~\eqref{eq2:thm:weak int log average thm 7} is proved from Eqs.~\eqref{eq4:thm:weak int log average thm 7} and~\eqref{eq:42}.

Next, we consider that $\mathcal{C}, \mathcal{D}_\tau$ are nonnegative Hermitian tensors. For any $\delta > 0$, we have following log-majorization relation:
\begin{eqnarray}
\prod\limits_{i=1}^k \left( \lambda_i (\mathcal{C}) + \epsilon_{\delta} \right) 
&\leq& \prod\limits_{i=1}^k \exp  \int_{\Omega} \log \left( \lambda_i(\mathcal{D}_\tau) + \delta\right) d \nu (\tau),
\end{eqnarray}
where $\delta >  \epsilon_{\delta} > 0$ and $k \in \{1,2,\cdots r \}$. Then, we can apply the previous case result about positive Hermitian tensors to positive Hermitian tensors $\mathcal{C} + \epsilon_{\delta} \mathcal{I}$ and $\mathcal{D}_\tau + \delta \mathcal{I}$, and get 
\begin{eqnarray}\label{eq:46}
\left\Vert f (\mathcal{C}) + \epsilon_{\delta}  \mathcal{I} \right\Vert_{(k)} 
&\leq& \exp \int_{\Omega} \log \left\Vert f (\mathcal{D}_{\tau}) + \delta  \mathcal{I} \right\Vert_{(k)} 
d \nu (\tau)
\end{eqnarray}
As $\delta \rightarrow 0$, Eq.~\eqref{eq:46} will give us Eq.~\eqref{eq2:thm:weak int log average thm 7} for nonnegative Hermitian tensors by the monotone
convergence theorem.  

\textbf{Eq.~\eqref{eq1:thm:weak int log average thm 7} $\Longleftarrow$ Eq.~\eqref{eq2:thm:weak int log average thm 7}}

We consider positive Hermitian tensors at first phase by assuming that $\mathcal{D}_{\tau}$ are 
positive Hermitian for all $\tau \in \Omega$. We may also assume that the tensor $\mathcal{C}$ is a positive Hermitian tensor. Since if this is a nonnegative Hermitian tensor, i.e., some $\lambda_i = 0$, we always have following inequality valid:
\begin{eqnarray}
\prod\limits_{i=1}^k \lambda_i (\mathcal{C}) \leq \prod\limits_{i=1}^k 
\exp \int_{\Omega} \log \lambda_i (\mathcal{D}_{\tau}) d \nu (\tau)
\end{eqnarray}

If we apply $f(x) = x^p$ for $p > 0$ and $\left\Vert \cdot \right\Vert_{(k)}$ as Ky Fan $k$-norm in Eq.~\eqref{eq2:thm:weak int log average thm 7}, we have 
\begin{eqnarray}\label{eq:50}
\log \sum\limits_{i=1}^k \lambda^p_i \left(\mathcal{C}\right) \leq \int_{\Omega} \log \sum\limits_{i=1}^k \lambda_i^p\left( \mathcal{D}_{\tau} \right) d \nu (\tau).
\end{eqnarray}
If we add $\log \frac{1}{k}$ and multiply $\frac{1}{p}$ at both sides of Eq.~\eqref{eq:50}, we have 
\begin{eqnarray}\label{eq:51}
\frac{1}{p}\log \left( \frac{1}{k}\sum\limits_{i=1}^k \lambda^p_i \left(\mathcal{C}\right) \right)\leq \int_{\Omega} \frac{1}{p} \log \left( \frac{1}{k}\sum\limits_{i=1}^k \lambda_i^p\left( \mathcal{D}_{\tau} \right) \right) d \nu (\tau).
\end{eqnarray}
From L'Hopital's Rule, if $p \rightarrow 0$, we have 
\begin{eqnarray}\label{eq:52}
\frac{1}{p}\log \left( \frac{1}{k}\sum\limits_{i=1}^k \lambda^p_i \left(\mathcal{C}\right) \right) \rightarrow \frac{1}{k} \sum\limits_{i=1}^k \log \lambda_i (\mathcal{C}),
\end{eqnarray}
and 
\begin{eqnarray}\label{eq:53}
\frac{1}{p}\log \left( \frac{1}{k}\sum\limits_{i=1}^k \lambda^p_i \left(\mathcal{D}_{\tau}\right) \right) \rightarrow \frac{1}{k} \sum\limits_{i=1}^k \log \lambda_i (\mathcal{D}_{\tau}),
\end{eqnarray}
where $\tau \in \Omega$. Appling Eqs.~\eqref{eq:52} and~\eqref{eq:53} into Eq.~\eqref{eq:51} and taking $p \rightarrow 0$, we have 
\begin{eqnarray}
\sum\limits_{i=1}^k \lambda_i (\mathcal{C}) \leq \int_{\Omega} \sum\limits_{i=1}^k 
 \log \lambda_i (\mathcal{D}_{\tau}) d \nu (\tau).
\end{eqnarray}
Therefore, Eq.~\eqref{eq1:thm:weak int log average thm 7} is true for positive Hermitian tensors. 

For nonnegative Hermitian tensors $\mathcal{D}_{\tau}$, since Eq.~\eqref{eq2:thm:weak int log average thm 7} is valid for $\mathcal{D}_{\tau} + \delta \mathcal{I}$ for any $\delta > 0$, we can apply the previous case result about positive Hermitian tensors to $\mathcal{D}_{\tau} + \delta \mathcal{I}$ and obtain
\begin{eqnarray}
\prod\limits_{i=1}^k \lambda_i (\mathcal{C}) \leq \prod\limits_{i=1}^k  \exp 
\int_{\Omega}   \log  \left( \lambda_i (\mathcal{D}_{\tau}) + \delta \right) d \nu (\tau),
\end{eqnarray}
where $k \in \{1,2,\cdots, r\}$. Eq.~\eqref{eq1:thm:weak int log average thm 7} is still true for nonnegative Hermitian tensors as $\delta \rightarrow 0$.

\textbf{Eq.~\eqref{eq1:thm:weak int log average thm 7} $\Longrightarrow$ Eq.~\eqref{eq3:thm:weak int log average thm 7}}

If $\mathcal{C}, \mathcal{D}_\tau$ are positive Hermitian tensors, and $\mathcal{D}_{\tau} \geq \delta \mathcal{I}$ for all $\tau \in \Omega$. From Eq.~\eqref{eq1:thm:weak int log average thm 7}, we have 
\begin{eqnarray}
\vec{\lambda} (\log \mathcal{C}) = \log \vec{\lambda}(\mathcal{C}) \prec_{w}
\int_{\Omega^r} \log \vec{\lambda}(\mathcal{D}_{\tau}) d \nu^r (\tau) = 
\int_{\Omega^r} \vec{\lambda}( \log \mathcal{D}_{\tau}) d \nu^r (\tau).
\end{eqnarray}
If we apply Theorem~\ref{thm:weak int average thm 4} to $\log \mathcal{C}$, $\log \mathcal{D}_{\tau}$ with function $f(x) = g(e^x)$, where $g$ is used in Eq.~\eqref{eq3:thm:weak int log average thm 7}, Eq.~\eqref{eq3:thm:weak int log average thm 7} is implied. 

If $\mathcal{C}, \mathcal{D}_\tau$ are nonnegative Hermitian tensors and any $\delta > 0$, we can find $\epsilon_{\delta} \in (0, \delta)$ to satisfy following:
\begin{eqnarray}\label{eq:45}
\prod\limits_{i=1}^k\left(\lambda_i(\mathcal{C}) + \epsilon_{\delta}\right) \leq 
\prod\limits_{i=1}^k \exp \int_{\Omega}   \log \left( \lambda_i(\mathcal{D}_{\tau}) + \delta  \right) d \nu (\tau).
\end{eqnarray}
Then, from positive Hermitian tensor case, we have 
\begin{eqnarray}\label{eq:45-1}
\left\Vert g( \mathcal{C} + \epsilon_{\delta} \mathcal{I} ) \right\Vert_{(k)}
\leq \int_{\Omega} \left\Vert   g( \mathcal{D}_{\tau} + \delta \mathcal{I} )    \right\Vert_{(k)}
d \nu (\tau).
\end{eqnarray}
Eq.~\eqref{eq3:thm:weak int log average thm 7} is obtained by taking $\delta \rightarrow 0$ in Eq.~\eqref{eq:45-1}. 

\textbf{Eq.~\eqref{eq1:thm:weak int log average thm 7} $\Longleftarrow$ Eq.~\eqref{eq3:thm:weak int log average thm 7}}

For $k \in \{1,2,\cdots, r \}$, if we apply $g(x) = \log (\delta + x )$, where $\delta >0$, and Ky Fan $k$-norm in Eq.~\eqref{eq3:thm:weak int log average thm 7}, we have 
\begin{eqnarray}
\sum\limits_{i=1}^k \log \left(\delta + \lambda_i \left(\mathcal{C} \right) \right)
\leq \sum\limits_{i=1}^k \int_{\Omega} \log \left( \delta + \lambda_{i}(\mathcal{D}_{\tau}) \right) d \nu (\tau).
\end{eqnarray}
Then, we have following as $\delta \rightarrow 0$:
\begin{eqnarray}
\sum\limits_{i=1}^k \log \lambda_i \left(\mathcal{C} \right) 
\leq \sum\limits_{i=1}^k \int_{\Omega} \log  \lambda_{i}(\mathcal{D}_{\tau}) d \nu (\tau).
\end{eqnarray}
Therefore, Eq.~\eqref{eq1:thm:weak int log average thm 7} ccan be derived from Eq.~\eqref{eq3:thm:weak int log average thm 7}. $\hfill \Box$

Next theorem will extend Theorem~\ref{thm:weak int log average thm 7} to non-weak version.

\begin{theorem}\label{thm:int log average thm 10}
Let $\mathcal{C}, \mathcal{D}_\tau$ be nonnegative Hermitian tensors with $\int_{\Omega} \left\Vert \mathcal{D}_{\tau}^{-p}\right\Vert_\rho d \nu (\tau) < \infty$ for any $p > 0$, $f: (0, \infty) \rightarrow [0,\infty)$ be a continous function such that the mapping $x \rightarrow \log f(e^x)$ is convex on $\mathbb{R}$, and $g: (0, \infty) \rightarrow [0,\infty)$ be a continous function such that the mapping $x \rightarrow g(e^x)$ is convex on $\mathbb{R}$ , then  we have following three equivalent statements:
\begin{eqnarray}\label{eq1:thm:int average thm 10}
\vec{\lambda}(\mathcal{C}) &\prec_{\log}& \exp  \int_{\Omega^r} \log \vec{\lambda}(\mathcal{D}_\tau) d\nu^r(\tau);
\end{eqnarray}
\begin{eqnarray}\label{eq2:thm:int average thm 10}
\left\Vert f(\mathcal{C}) \right\Vert_{(k)} &\leq &
\exp \int_{\Omega} \log \left\Vert f(\mathcal{D}_{\tau}) \right\Vert_{(k)}  d\nu(\tau);
\end{eqnarray}
\begin{eqnarray}\label{eq3:thm:int average thm 10}
\left\Vert g(\mathcal{C}) \right\Vert_{(k)} &\leq &
\int_{\Omega} \left\Vert g(\mathcal{D}_{\tau}) \right\Vert_{(k)}  d\nu(\tau).
\end{eqnarray}
\end{theorem}
\textbf{Proof:}

The proof plan is similar to the proof in Theorem~\ref{thm:weak int log average thm 7}. 

\textbf{Eq.~\eqref{eq1:thm:int average thm 10} $\Longrightarrow$ Eq.~\eqref{eq2:thm:int average thm 10}}

First, we assume that $\mathcal{C}, \mathcal{D}_\tau$ are postiive Hermitian tensors with $\mathcal{D}_\tau \geq \delta \mathcal{I}$ for all $\tau \in \Omega$. The corresponding part of the proof in Theorem~\ref{thm:weak int log average thm 7} about positive Hermitian tensors $\mathcal{C}, \mathcal{D}_\tau$ can be applied here. 

For case that $\mathcal{C}, \mathcal{D}_\tau$ are nonnegative Hermitian tensors, we have 
\begin{eqnarray}
\prod\limits_{i=1}^k \lambda_i (\mathcal{C}) \leq \prod\limits_{i=1}^k 
\exp \int_{\Omega} \log \left( \lambda_i (\mathcal{D}_{\tau}) + \delta_n \right)d \nu (\tau), 
\end{eqnarray}
where $\delta_n > 0$ and $\delta_n \rightarrow 0$. Because $\int_{\Omega^r} \log \left( \vec{\lambda} (\mathcal{D}_{\tau}) + \delta_n \right) d \nu^r (\tau) \rightarrow \int_{\Omega^r} \log \vec{\lambda} (\mathcal{D}_\tau)  d \nu^r (\tau) $ as $n \rightarrow \infty$, from Lemma~\ref{lma:Lemma 12 Gen Log Hiai}, we can find $\mathbf{a}^{(n)}$ with $n \geq n_0$ such that $a^{(n)}_1 \geq \cdots \geq a^{(n)}_r > 0$, $\mathbf{a}^{(n)} \rightarrow \vec{\lambda}(\mathcal{C})$ and $ \mathbf{a}^{(n)} \prec_{\log}  \exp \int_{\Omega^r} \log \vec{\lambda} \left(\mathcal{D}_{\tau} + \delta_n \mathcal{I} \right) d \nu^r (\tau)$

Selecting $\mathcal{C}^{(n)}$ with $\vec{\lambda} ( \mathcal{C}^{(n)})  = \mathbf{a}^{(n)} $ and applying positive Hermitian tensors case to $\mathcal{C}^{(n)}$ and $\mathcal{D}_{\tau} + \delta_n \mathcal{I}$, we obtain
\begin{eqnarray}\label{eq:74}
\left\Vert f (\mathcal{C}^{(n)}) \right\Vert_{(k)} \leq \exp \int_{\Omega} \log \left\Vert f (\mathcal{D}_{\tau} + \delta_n \mathcal{I}) \right\Vert_{(k)} d \nu (\tau)
\end{eqnarray}
where $n \geq n_0$.

There are two situations for the function $f$ near $0$: $f(0^{+}) < \infty$ and $f(0^{+}) = \infty$. For the case with $f(0^{+}) < \infty$, we have 
\begin{eqnarray}\label{eq:75-1}
\left\Vert f (\mathcal{C}^{(n)} )\right\Vert_{(k)} = \rho( f (\mathbf{a}^{(n)}))
\rightarrow \rho (f ( \vec{\lambda}(\mathcal{C}))) = \left\Vert f (\mathcal{C})\right\Vert_{(k)}, 
\end{eqnarray}
and
\begin{eqnarray}\label{eq:75-2}
\left\Vert f (\mathcal{D}_{\tau} + \delta_n \mathcal{I} )\right\Vert_{(k)} 
\rightarrow \left\Vert f (\mathcal{D}_{\tau})\right\Vert_{(k)}, 
\end{eqnarray}
where $\tau \in \Omega$ and $n \rightarrow \infty$. From Fatou–Lebesgue theorem, we then have 
\begin{eqnarray}\label{eq:76}
\limsup\limits_{n \rightarrow \infty} \int_{\Omega} \log \left\Vert f (\mathcal{D}_{\tau} + \delta_n \mathcal{I} )\right\Vert_{(k)} d \nu (\tau) \leq \int_{\Omega} \log \left\Vert f(\mathcal{D}_{\tau}) \right\Vert_{(k)}.
\end{eqnarray}
By taking $n \rightarrow \infty$ in Eq.~\eqref{eq:74} and using Eqs.~\eqref{eq:75-1},~\eqref{eq:75-2},~\eqref{eq:76}, we have Eq.~\eqref{eq2:thm:int average thm 10} for case that $f(0^{+}) < \infty$.

For the case with $f(0^{+}) = \infty$, we assume that $\int_{\Omega} \log \left\Vert f (\mathcal{D}_{\tau}) \right\Vert_{(k)} d \nu (\tau) < \infty$ (since the inequality in Eq.~\eqref{eq2:thm:int average thm 10} is always true for $\int_{\Omega} \log \left\Vert f (\mathcal{D}_{\tau}) \right\Vert_{(k)} d \nu (\tau) = \infty$). Since $f$ is decreasing on $(0, \epsilon)$ for some $\epsilon > 0$. We claim that the following relation is valid: there are two constants $a, b > 0$ such that 
\begin{eqnarray}\label{eq:77}
a \leq \left\Vert f (\mathcal{D}_{\tau} + \delta_n \mathcal{I}) \right\Vert_{(k)}
\leq \left\Vert f (\mathcal{D}_{\tau}) \right\Vert_{(k)} + b,
\end{eqnarray}
for all $\tau \in \Omega$ and $n \geq n_0$. If Eq.~\eqref{eq:77} is valid and $\int_{\Omega} \log \left\Vert f (\mathcal{D}_{\tau}) \right\Vert_{(k)} d \nu (\tau) < \infty$, from Lebesgue's dominated convergence theorem, we also have Eq.~\eqref{eq2:thm:int average thm 10} for case that $f(0^{+}) = \infty$ by taking $n \rightarrow \infty$ in Eq.~\eqref{eq:74}. 

Now, we want to prove the claim shown by Eq.~\eqref{eq:77}. By the uniform boundedness of tensors $\mathcal{D}_{\tau}$, there is a constant $\chi >0$ such that 
\begin{eqnarray}
0 < \mathcal{D}_{\tau} + \delta_n \mathcal{I} \leq \chi \mathcal{I},
\end{eqnarray}
where $\tau \in \Omega$ and $ n \geq n_0$. 
%
From Eq.~\eqref{eq:Hermitian Eigen Decom}, we have 
\begin{eqnarray}
f(\mathcal{D}_{\tau} + \delta_n \mathcal{I}) &=&  \sum\limits_{i', \mbox{s.t. $\lambda_{i'}(\mathcal{D}_{\tau})+  \delta_n < \epsilon$}} f(\lambda_{i'}(\mathcal{D}_{\tau}) + \delta_n ) \mathcal{U}_{i'} \otimes \mathcal{U}^{H}_{i'} + \nonumber \\
&  & \sum\limits_{j', \mbox{s.t. $\lambda_{j'}(\mathcal{D}_{\tau}) +  \delta_n 
\geq \epsilon$}} f(\lambda_{j'}(\mathcal{D}_{\tau}) + \delta_n ) \mathcal{U}_{j'} \otimes \mathcal{U}^{H}_{j'}  \nonumber \\
&\leq &  \sum\limits_{i', \mbox{s.t. $\lambda_{i'}(\mathcal{D}_{\tau}) +  \delta_n < \epsilon$}} f(\lambda_{i'}(\mathcal{D}_{\tau}) ) \mathcal{U}_{i'} \otimes \mathcal{U}^{H}_{i'} + \nonumber \\
&  & \sum\limits_{j', \mbox{s.t. $\lambda_{j'}(\mathcal{D}_{\tau})  +  \delta_n 
\geq \epsilon$}} f(\lambda_{j'}(\mathcal{D}_{\tau}) + \delta_n ) \mathcal{U}_{j'} \otimes \mathcal{U}^{H}_{j'}  \nonumber \\
&\leq & f(\mathcal{D}_{\tau}) + \sum\limits_{j', \mbox{s.t. $\lambda_{j'}(\mathcal{D}_{\tau}) +  \delta_n \geq \epsilon$}} f(\lambda_{j'}(\mathcal{D}_{\tau}) + \delta_n ) \mathcal{U}_{j'} \otimes \mathcal{U}^{H}_{j'}.
\end{eqnarray}
Therefore, the claim in Eq.~\eqref{eq:77} follows by the triangle inequality for $\left\Vert \cdot \right\Vert_{(k)}$ and $f(\lambda_{j'}(\mathcal{D}_{\tau}) + \delta_n )  < \infty$ for $\lambda_{j'}(\mathcal{D}_{\tau}) +  \delta_n \geq \epsilon$. 

\textbf{Eq.~\eqref{eq1:thm:int average thm 10} $\Longleftarrow$ Eq.~\eqref{eq2:thm:int average thm 10}}

The weak majorization relation 
\begin{eqnarray}\label{eq:82}
\prod\limits_{i=1}^{k} \lambda_i (\mathcal{C}) \leq \prod\limits_{i=1}^{k} \exp \int_{\Omega} \log \lambda_i (\mathcal{D}_\tau) d \nu (\tau),
\end{eqnarray}
is valid for $k < r$ from Eq.~\eqref{eq1:thm:weak int log average thm 7} $\Longrightarrow$ Eq.~\eqref{eq2:thm:weak int log average thm 7} in Theorem~\ref{thm:weak int log average thm 7}.   We wish to prove that Eq.~\eqref{eq:82} becomes equal for $k = r$. It is equivalent to prove that 
\begin{eqnarray}\label{eq:83}
\log \det\nolimits_H (\mathcal{C}) \geq   \int_{\Omega} \log \det\nolimits_H (\mathcal{D}_{\tau}) d \nu (\tau),
\end{eqnarray}
where $\det\nolimits_H ( \cdot )$ is the \emph{Hermitian determinant}. We can assume that $ \int_{\Omega} \log \det\nolimits_H (\mathcal{D}_{\tau}) d \nu (\tau) \geq - \infty$ since Eq.~\eqref{eq:83} is true for  $ \int_{\Omega} \log \det\nolimits_H (\mathcal{D}_{\tau}) d \nu (\tau) = - \infty$. Then, $\mathcal{D}_\tau$ are positive Hermitian tensors. 

If we scale tensors $\mathcal{C}, \mathcal{D}_{\tau}$ as $a \mathcal{C}, a\mathcal{D}_{\tau}$ by some $a >0$, we can assume $\mathcal{D}_{\tau} \leq \mathcal{I}$ and $\lambda_i(\mathcal{D}_\tau) \leq 1$ for all $ \tau \in \Omega$ and $ i \in \{1,2,\cdots, r\}$. Then for any $p >0$, we have 
\begin{eqnarray}
\frac{1}{r} \left\Vert \mathcal{D}_{\tau}^{-p} \right\Vert_1 \leq \lambda^{-p}_r (\mathcal{D}_{\tau} ) \leq ( \det\nolimits_H  (\mathcal{D}_{\tau})  )^{-p},
\end{eqnarray}
where $\left\Vert \cdot \right\Vert_1$ represents the tensor trace norm, and 
\begin{eqnarray}\label{eq:85}
\frac{1}{p} \log \left( \frac{\left\Vert \mathcal{D}^{-p}_\tau \right\Vert_1 }{r}\right)
\leq - \log \det\nolimits_H  (\mathcal{D}_{\tau}). 
\end{eqnarray}
If we use tensor trace norm as unitarily invariant tensor norm and $f(x) = x^{-p}$ for any $p > 0$ in Eq.~\eqref{eq2:thm:int average thm 10}, we obtain
\begin{eqnarray}\label{eq:86}
\log \left\Vert \mathcal{C}^{-p} \right\Vert_1 \leq \int_{\Omega} \log \left\Vert \mathcal{D}^{-p}_{\tau} \right\Vert_1 d \nu(\tau).
\end{eqnarray}
By adding $\log \frac{1}{r}$ and multiplying $\frac{1}{p}$ for both sides of Eq.~\eqref{eq:86}, we have 
\begin{eqnarray}\label{eq:87}
\frac{1}{p}\log \left( \frac{\left\Vert \mathcal{C}^{-p} \right\Vert_1 }{r} \right)
\leq \int_{\Omega}\frac{1}{p} \log \left(  \frac{\left\Vert \mathcal{D}_\tau^{-p} \right\Vert_1 }{r}           \right) d \nu (\tau)
\end{eqnarray}
Similar to Eqs.~\eqref{eq:52} and~\eqref{eq:53}, we have following two relations as $p \rightarrow 0$:
\begin{eqnarray}\label{eq:88}
\frac{1}{p}\log \left( \frac{\left\Vert \mathcal{C}^{-p} \right\Vert_1 }{r} \right) \rightarrow \frac{-  \log \det\nolimits_H (\mathcal{C}) }{r},
\end{eqnarray}
and 
\begin{eqnarray}\label{eq:89}
\frac{1}{p}\log \left( \frac{\left\Vert \mathcal{D}_{\tau}^{-p} \right\Vert_1}{r}  \right) \rightarrow \frac{- \log \det\nolimits_H (\mathcal{D}_\tau) }{r} .
\end{eqnarray}
From Eq.~\eqref{eq:85} and Lebesgue's dominated convergence theorem, we have 
\begin{eqnarray}\label{eq:90}
\lim\limits_{p  \rightarrow 0} \int_{\Omega}\frac{1}{p}\log \left( \frac{\left\Vert \mathcal{D}_{\tau}^{-p} \right\Vert_1}{r}  \right) d \nu (\tau)= \frac{-1}{r} \int_{\Omega} \log \det\nolimits_{H}(\mathcal{D}_{\tau})     \nu (\tau) 
\end{eqnarray}
Finally, we have Eq.~\eqref{eq:83} from Eqs.~\eqref{eq:87} and~\eqref{eq:90}.

\textbf{Eq.~\eqref{eq1:thm:int average thm 10} $\Longrightarrow$ Eq.~\eqref{eq3:thm:int average thm 10}}

First, we assume that $\mathcal{C}, \mathcal{D}_\tau$ are positive Hermitian tensors and $\mathcal{D}_{\tau} \geq \delta \mathcal{I}$ for $\tau \in \Omega$. From Eq.~\eqref{eq1:thm:int average thm 10}, we can apply Theorem~\ref{thm:weak int average thm 5} to $\log \mathcal{C}, \log \mathcal{D}_\tau$ and $f(x) = g(e^x)$ to obtain Eq.~\eqref{eq3:thm:int average thm 10}. 

For $\mathcal{C}, \mathcal{D}_\tau$ are nonnegative Hermitian tensors, we can choose  $\mathbf{a}^{(n)}$ and corresponding  $\mathcal{C}^{(n)}$ for $n \geq n_0$ given $\delta_n \rightarrow 0$ with $\delta_n > 0$ as the proof in Eq.~\eqref{eq1:thm:int average thm 10} $\Longrightarrow$ Eq.~\eqref{eq2:thm:int average thm 10}. Since tensors $\mathcal{C}^{(n)}, \mathcal{D}_\tau + \delta_n \mathcal{I}$ are postive Hermitian tensors, we then have 
\begin{eqnarray}\label{eq:92}
\left \Vert g ( \mathcal{C}^{(n)}  ) \right\Vert_{(k)} \leq \int_{\Omega} \left\Vert g ( \mathcal{D}_\tau + \delta_n \mathcal{I} ) \right\Vert_{(k)} d \nu (\tau).
\end{eqnarray}
If $g(0^+) < \infty$, Eq.~\eqref{eq3:thm:int average thm 10} is obtained from Eq.~\eqref{eq:92} by taking $n \rightarrow \infty$. On the other hand, if $g(0^+) = \infty$, we can apply the argument similar to the portion about $f(0^+) = \infty$ in the proof for Eq.~\eqref{eq1:thm:int average thm 10} $\Longrightarrow$ Eq.~\eqref{eq2:thm:int average thm 10} to get $a, b > 0$ such that 
\begin{eqnarray}\label{eq:92 infty}
a \leq \left\Vert g ( \mathcal{D}_\tau + \delta_n \mathcal{I} ) \right\Vert_\rho \leq   \left\Vert  g ( \mathcal{D}_\tau ) \right\Vert_\rho + b,
\end{eqnarray}
for all $\tau \in \Omega$ and $n \geq n_0$. Since the case that $\int_{\Omega} \left\Vert  g ( \mathcal{D}_\tau ) \right\Vert_{(k)} d \nu (\tau) = \infty$ will have Eq.~\eqref{eq3:thm:int average thm 10}, we only consider the case that $\int_{\Omega} \left\Vert  g ( \mathcal{D}_\tau ) \right\Vert_{(k)} d \nu (\tau) < \infty$. Then, we have Eq.~\eqref{eq3:thm:int average thm 10} from Eqs.~\eqref{eq:92},~\eqref{eq:92 infty} and Lebesgue's dominated convergence theorem.

\textbf{Eq.~\eqref{eq1:thm:int average thm 10} $\Longleftarrow$ Eq.~\eqref{eq3:thm:int average thm 10}}

The weak majorization relation
\begin{eqnarray}\label{eq:94}
\sum\limits_{i=1}^{k} \log \lambda_i (\mathcal{C}) \leq 
\sum\limits_{i=1}^{k} \int_{\Omega} \log \lambda_i (\mathcal{D}_\tau) d \nu (\tau)
\end{eqnarray}
is true from the implication from Eq.~\eqref{eq3:thm:weak int log average thm 7} to Eq.~\eqref{eq1:thm:weak int log average thm 7} in Theorem~\ref{thm:weak int log average thm 7}. We have to show that this relation becomes identity for $k=r$. If we apply $\left\Vert \cdot  \right\Vert_{\rho} = \left\Vert \cdot \right\Vert_1$ and $g(x) = x^{-p}$ for any $p > 0$ in Eq.~\eqref{eq3:thm:int average thm 10}, we have 
\begin{eqnarray}\label{eq:95}
\frac{1}{p} \log \left( \frac{ \left\Vert \mathcal{C}^{-p}\right\Vert_1  }{r} \right)
\leq \frac{1}{p} \log \int_{\Omega} \frac{\left\Vert \mathcal{D}_\tau^{-p}\right\Vert_1}{r}   d \nu(\tau).
\end{eqnarray}
Then, we will get 
\begin{eqnarray}\label{eq:96}
\frac{- \log \det\nolimits_H (\mathcal{C})}{r} &=& \lim\limits_{p \rightarrow 0} \frac{1}{p} \log \left( \frac{ \left\Vert \mathcal{C}^{-p}\right\Vert_1  }{r} \right) \nonumber \\
& \leq & \lim\limits_{p \rightarrow 0} \frac{1}{p} \log \left( \int_{\Omega} \frac{\left\Vert \mathcal{D}_\tau^{-p}\right\Vert_1}{r}   d \nu(\tau)  \right)=_1 \frac{ - \int_{\Omega} \log \det\nolimits_H (\mathcal{D}_{\tau}) d \nu (\tau)  }{r},
\end{eqnarray}
which will proves the identity for Eq.~\eqref{eq:94} when $k = r$. The equality in $=_1$ will be proved by the following Lemma~\ref{lma:15}.
$\hfill \Box$

\begin{lemma}\label{lma:15}
Let $\mathcal{D}_\tau$ be nonnegative Hermitian tensors with $\int_{\Omega} \left\Vert \mathcal{D}_{\tau}^{-p}\right\Vert_1 d \nu (\tau) < \infty$ for any $p > 0$, then we have
\begin{eqnarray}\label{eq1:lma:15}
\lim\limits_{p \rightarrow 0} \left( \frac{1}{p} \log \int_{\Omega} \frac{ \left\Vert \mathcal{D}_\tau^{-p}\right\Vert_1 }{r} d \nu (\tau)\right) &=& -\frac{  \int_{\Omega} \log \det\nolimits_{H}( \mathcal{D}_{\tau} ) d \nu (\tau) }{r} 
\end{eqnarray}
\end{lemma}
\textbf{Proof:}
Because $\int_{\Omega} \left\Vert \mathcal{D}_{\tau}^{-p}\right\Vert_\rho d \nu (\tau) < \infty$, we have that $\mathcal{D}_{\tau}$ are positive Hermitian tensors for $\tau$ almost everywhere in $\Omega$. Then, we have 
\begin{eqnarray}
\lim\limits_{p \rightarrow 0} \left( \frac{1}{p}\log \int_{\Omega} \frac{ \left\Vert \mathcal{D}_\tau^{-p}\right\Vert_1}{r} d \nu (\tau) \right) &=_1& \lim\limits_{p \rightarrow 0}\frac{   \int_{\Omega} \frac{ - \sum\limits_{i=1}^r \log \lambda_i(\mathcal{D}_\tau)   }{r} d \nu (\tau)   }{  \int_{\Omega} \frac{ \left\Vert \mathcal{D}_\tau^{-p}\right\Vert_1}{r} d \nu (\tau)      } \nonumber \\
&=& \frac{-1}{r} \int_{\Omega} \sum\limits_{i=1}^r \log \lambda_i(\mathcal{D}_\tau)  d \nu (\tau)  
\nonumber \\
&=_2& \frac{-1}{r} \int_{\Omega} \log \det\nolimits_H ( \mathcal{D}_\tau ) d \nu (\tau), 
\end{eqnarray}
where $=_1$ is from L'Hopital's rule, and $=_2$ is obtained from $\det\nolimits_H$ definition.
$\hfill \Box$

\subsection{Multivaraite Tensor Norm Inequalities}\label{sec:Multivaraite Tensor Norm Inequalities}



In this section, we will apply derived majorization inequalities for tensors to multivaraite tensor norm inequalities which will be used to derive tensor expander bounds. We will begin to present a Lie-Trotter product formula for tensors. 
\begin{lemma}\label{lma: Lie product formula for tensors}
Let $m \in \mathbb{N}$ and $(\mathcal{L}_k)_{k=1}^{m}$ be a finite sequence of bounded tensors with dimensions $\mathcal{L}_k \in  \mathbb{C}^{I_1 \times \cdots \times I_M\times I_1 \times \cdots \times I_M}$, then we have
\begin{eqnarray}
\lim_{n \rightarrow \infty} \left(  \prod_{k=1}^{m} \exp(\frac{  \mathcal{L}_k}{n})\right)^{n}
&=& \exp \left( \sum_{k=1}^{m}  \mathcal{L}_k \right)
\end{eqnarray}
\end{lemma}
\textbf{Proof:}

We will prove the case for $m=2$, and the general value of $m$ can be obtained by mathematical induction. 
Let $\mathcal{L}_1, \mathcal{L}_2$ be bounded tensors act on some Hilbert space. Define $\mathcal{C} \define \exp( (\mathcal{L}_1 + \mathcal{L}_2)/n) $, and $\mathcal{D} \define \exp(\mathcal{L}_1/n) \star_M \exp(\mathcal{L}_2/n)$. Note we have following estimates for the norm of tensors $\mathcal{C}, \mathcal{D}$: 
\begin{eqnarray}\label{eq0: lma: Lie product formula for tensors}
\left\Vert \mathcal{C} \right\Vert, \left\Vert \mathcal{D} \right\Vert \leq \exp \left( \frac{\left\Vert \mathcal{L}_1 \right\Vert + \left\Vert \mathcal{L}_2 \right\Vert  }{n} \right) =  \left[ \exp \left(  \left\Vert \mathcal{L}_1 \right\Vert + \left\Vert \mathcal{L}_2 \right\Vert  \right) \right]^{1/n}.
\end{eqnarray}

From the Cauchy-Product formula, the tensor $\mathcal{D}$ can be expressed as:
\begin{eqnarray}
\mathcal{D} &=& \exp(\mathcal{L}_1/n) \star_M \exp(\mathcal{L}_2/n) = \sum_{i = 0}^{\infty} \frac{( \mathcal{L}_1/n)^i}{i !} \star_M \sum_{j = 0}^{\infty} \frac{( \mathcal{L}_2/n)^j}{j !} \nonumber\\
&=& \sum_{l = 0}^{\infty} n^{-l} \sum_{i=0}^l \frac{\mathcal{L}_1^i}{i!} \star_M \frac{\mathcal{L}_2^{l-i}}{(l - i)!},
\end{eqnarray}
then we can bound the norm of $\mathcal{C} - \mathcal{D}$ as 
\begin{eqnarray}\label{eq1: lma: Lie product formula for tensors}
\left\Vert \mathcal{C} - \mathcal{D} \right\Vert &=& \left\Vert \sum_{i=0}^{\infty} \frac{\left( [ \mathcal{L}_1 + \mathcal{L}_2]/n \right)^i}{i! }
 - \sum_{l = 0}^{\infty} n^{-l} \sum_{i=0}^l \frac{\mathcal{L}_1^i}{i!} \star_M \frac{\mathcal{L}_2^{l-i}}{(l - i)!} \right\Vert \nonumber \\
&=& \left\Vert \sum_{i=2}^{\infty} k^{-i} \frac{\left( [ \mathcal{L}_1 + \mathcal{L}_2] \right)^i}{i! }
 - \sum_{m = l}^{\infty} n^{-l} \sum_{i=0}^l \frac{\mathcal{L}_1^i}{i!} \star_M \frac{\mathcal{L}_2^{l-i}}{(l - i)!} \right\Vert \nonumber \\
& \leq & \frac{1}{k^2}\left[ \exp( \left\Vert \mathcal{L}_1 \right\Vert + \left\Vert \mathcal{L}_2 \right\Vert ) + \sum_{l = 2}^{\infty} n^{-l} \sum_{i=0}^l \frac{\left\Vert \mathcal{L}_1 \right\Vert^i}{i!} \cdot \frac{\left\Vert \mathcal{L}_2 \right\Vert^{l-i}}{(l - i)!} \right] \nonumber \\
& = & \frac{1}{n^2}\left[ \exp \left( \left\Vert \mathcal{L}_1 \right\Vert + \left\Vert \mathcal{L}_2 \right\Vert \right) + \sum_{l = 2}^{\infty} n^{-l} \frac{(  \left\Vert \mathcal{L}_1 \right\Vert + \left\Vert \mathcal{L}_2 \right\Vert )^l}{l!} \right] \nonumber \\
& \leq & \frac{2  \exp \left( \left\Vert \mathcal{L}_1 \right\Vert + \left\Vert \mathcal{L}_2 \right\Vert \right) }{n^2}.
\end{eqnarray}

For the difference between the higher power of $\mathcal{C}$ and $\mathcal{D}$, we can bound them as 
\begin{eqnarray}
\left\Vert \mathcal{C}^n - \mathcal{D}^n \right\Vert &=& \left\Vert \sum_{l=0}^{n-1} \mathcal{C}^m (\mathcal{C} - \mathcal{D})\mathcal{C}^{n-l-1} \right\Vert \nonumber \\
& \leq_1 &  \exp ( \left\Vert \mathcal{L}_1 \right\Vert +  \left\Vert \mathcal{L}_2 \right\Vert) \cdot n \cdot \left\Vert \mathcal{L}_1 - \mathcal{L}_2 \right\Vert,
\end{eqnarray}
where the inequality $\leq_1$ uses the following fact 
\begin{eqnarray}
\left\Vert \mathcal{C} \right\Vert^{l} \left\Vert \mathcal{D} \right\Vert^{n - l - 1} \leq \exp \left( \left\Vert \mathcal{L}_1 \right\Vert +  \left\Vert \mathcal{L}_2 \right\Vert \right)^{\frac{n-1}{n}} \leq 
 \exp\left( \left\Vert \mathcal{L}_1 \right\Vert +  \left\Vert \mathcal{L}_2 \right\Vert \right), 
\end{eqnarray}
based on Eq.~\eqref{eq0: lma: Lie product formula for tensors}. By combining with Eq.~\eqref{eq1: lma: Lie product formula for tensors}, we have the following bound
\begin{eqnarray}
\left\Vert \mathcal{C}^n - \mathcal{D}^n \right\Vert &\leq& \frac{2 \exp \left( 2  \left\Vert \mathcal{L}_1 \right\Vert  +  2  \left\Vert \mathcal{L}_2 \right\Vert \right)}{n}.
\end{eqnarray}
Then this lemma is proved when $n$ goes to infity. $\hfill \Box$

\begin{theorem}\label{thm:Multivaraite Tensor Norm Inequalities}
Let $\mathcal{C}_i \in \mathbb{C}^{I_1 \times \cdots \times I_N \times I_1 \times \cdots \times I_N}$ be positive Hermitian tensors for $1 \leq i \leq n$ with Hermitian rank $r$, $\left\Vert \cdot \right\Vert_{(k)}$ be Ky Fan $k$-norm. For any continous function $f:(0, \infty) \rightarrow [0, \infty)$ such that $x \rightarrow \log f(e^x)$ is convex on $\mathbb{R}$, we have 
\begin{eqnarray}\label{eq1:thm:Multivaraite Tensor Norm Inequalities}
\left\Vert  f \left( \exp \left( \sum\limits_{i=1}^n \log \mathcal{C}_i\right)   \right)  \right\Vert_{(k)} &\leq& \exp \int_{- \infty}^{\infty} \log \left\Vert f \left( \left\vert \prod\limits_{i=1}^{n}  \mathcal{C}_i^{1 + \iota t} \right\vert\right)\right\Vert_{(k)} \beta_0(t) dt ,
\end{eqnarray}
where $\beta_0(t) = \frac{\pi}{2 (\cosh (\pi t) + 1)}$.

For any continous function $g(0, \infty) \rightarrow [0, \infty)$ such that $x \rightarrow g (e^x)$ is convex on $\mathbb{R}$, we have 
\begin{eqnarray}\label{eq2:thm:Multivaraite Tensor Norm Inequalities}
\left\Vert  g \left( \exp \left( \sum\limits_{i=1}^n \log \mathcal{C}_i\right)   \right)  \right\Vert_{(k)} &\leq& \int_{- \infty}^{\infty} \left\Vert g \left( \left\vert \prod\limits_{i=1}^{n}  \mathcal{C}_i^{1 + \iota t} \right\vert\right)\right\Vert_{(k)} \beta_0(t) dt.
\end{eqnarray}
\end{theorem}
\textbf{Proof:}
From Hirschman interpolation theorem~\cite{sutter2017multivariate} and $\theta \in [0, 1]$ and setting $\sqrt{-1}$ as $\iota$, we have 
\begin{eqnarray}\label{eq1:Hirschman interpolation}
\log \left\vert h(\theta) \right\vert \leq \int_{- \infty}^{\infty} \log \left\vert h(\iota t) \right\vert^{1 - \theta} \beta_{1 - \theta}(t) dt+ \int_{- \infty}^{\infty} \log \left\vert h(1 +  \iota t) \right\vert^{\theta} \beta_{\theta}(t) dt, 
\end{eqnarray}
where $h(z)$ be uniformly bounded on $S \define \{ z \in \mathbb{C}: 0 \leq \Re(z) \leq 1  \}$ and holomorphic on $S$. The term $ \beta_{\theta}(t) $ is defined as :
\begin{eqnarray}\label{eq:beta theta t def}
\beta_{\theta}(t) \define \frac{ \sin (\pi \theta)}{ 2 \theta (\cos(\pi t) + \cos (\pi \theta))  }.  
\end{eqnarray}
Let $H(z)$ be a uniformly bounded holomorphic function with values in $\mathbb{C}^{I_1 \times \cdots \times I_N \times I_1 \times \cdots \times I_N}$. Fix some $\theta \in [0, 1]$ and let $\mathcal{U}, \mathcal{V} \in \mathbb{C}^{I_1 \times \cdots \times I_N \times 1}$ be normalized tensors such that $\langle \mathcal{U}, \mathcal{H}(\theta) \star_N \mathcal{V} \rangle = \left\Vert H(\theta) \right\Vert$. If we define $h(z)$ as $h(z) \define \langle \mathcal{U}, \mathcal{H}(z) \star_N \mathcal{V} \rangle $, we have following bound: $\left\vert h(z) \right\vert \leq \left\Vert H(z) \right\Vert $ for all $z \in S$. From Hirschman interpolation theorem, we then have following interpolation theorem for tensor-valued function: 
\begin{eqnarray}\label{eq2:Hirschman interpolation}
\log \left\Vert H(\theta) \right\Vert \leq \int_{- \infty}^{\infty} \log \left\Vert H(\iota t) \right\Vert^{1 - \theta} \beta_{1 - \theta}(t) d t + \int_{- \infty}^{\infty} \log \left\Vert H(1 +  \iota t) \right\Vert^{\theta}  \beta_{\theta}(t) d t .  
\end{eqnarray}

Let $H(z) = \prod\limits_{i=1}^{n} \mathcal{C}^z_i$. Then the first term in the R.H.S. of Eq.~\eqref{eq2:Hirschman interpolation} is zero since $H(\iota t)$ is a product of unitary tensors. Then we have 
\begin{eqnarray}\label{eq:109}
\log \left\Vert \left\vert \prod\limits_{i=1}^{n} \mathcal{C}_i^{\theta} \right\vert^{\frac{1}{\theta}} \right\Vert \leq \int_{- \infty}^{\infty} \log \left\Vert    \prod\limits_{i=1}^{n} \mathcal{C}_i^{1 + \iota t}         \right\Vert  \beta_{\theta}(t) dt.  
\end{eqnarray}

From Lemma~\ref{lma:antisymmetric tensor product properties}, we have following relations:
\begin{eqnarray}\label{eq:111-1}
\left\vert \prod\limits_{i=1}^n \left( \wedge^k \mathcal{C}_i \right)^{\theta}\right\vert^{\frac{1}{\theta}} = \wedge^k \left\vert  \prod\limits_{i=1}^n \mathcal{C}^{\theta}_i  \right\vert^{\frac{1}{\theta}},
\end{eqnarray}
and
\begin{eqnarray}\label{eq:111-2}
\left\vert \prod\limits_{i=1}^n \left( \wedge^k \mathcal{C}_i \right)^{1 + \iota t} \right\vert = \wedge^k \left\vert  \prod\limits_{i=1}^n \mathcal{C}^{1 + \iota t}_i  \right\vert.
\end{eqnarray}
If Eq.~\eqref{eq:109} is applied to $\wedge^k \mathcal{C}_i$ for $1 \leq k \leq r$, we have following log-majorization relation from Eqs.~\eqref{eq:111-1} and~\eqref{eq:111-2}:
\begin{eqnarray}\label{eq:112}
\log \vec{\lambda} \left(  \left\vert \prod\limits_{i=1}^{n} \mathcal{C}_i^{\theta} \right\vert^{\frac{1}{\theta}}   \right) \prec \int_{- \infty}^{\infty} \log \vec{\lambda} \left\vert \prod\limits_{i=1}^{n} \mathcal{C}_i^{1 + \iota  t } \right\vert^{\frac{1}{\theta}}  \beta_{\theta}(t) d t .
\end{eqnarray}
Moreover, we have the equality condition in Eq.~\eqref{eq:112} for $k = r$ due to following identies:
\begin{eqnarray}\label{eq:113}
\det\nolimits_H \left\vert \prod\limits_{i=1}^n \mathcal{C}_i^{\theta} \right\vert^{\frac{1}{\theta}}
= \det\nolimits_H \left\vert \prod\limits_{i=1}^n \mathcal{C}_i^{1 + \iota t } \right\vert= \prod\limits_{i=1}^n \det\nolimits_H \mathcal{C}_i. 
\end{eqnarray}

At this stage, we are ready to apply Theorem~\ref{thm:int log average thm 10}  for the log-majorization provided by Eq.~\eqref{eq:112} to get following facts:
\begin{eqnarray}\label{eq:114}
\left\Vert  f \left( \left\vert \prod\limits_{i=1}^{n} \mathcal{C}_i^{\theta} \right\vert^{\frac{1}{\theta}}  \right)  \right\Vert_{(k)} &\leq& \exp \int_{- \infty}^{\infty} \log \left\Vert f \left( \left\vert \prod\limits_{i=1}^{n}  \mathcal{C}_i^{1 + \iota t} \right\vert\right)\right\Vert_{(k)}  \beta_{\theta}(t) d t  ,
\end{eqnarray}
and
\begin{eqnarray}\label{eq:115}
\left\Vert  g \left( \left\vert \prod\limits_{i=1}^{n} \mathcal{C}_i^{\theta} \right\vert^{\frac{1}{\theta}}  \right)  \right\Vert_{(k)} &\leq& \int_{- \infty}^{\infty} \left\Vert g \left( \left\vert \prod\limits_{i=1}^{n}  \mathcal{C}_i^{1 + \iota t} \right\vert\right)\right\Vert_{(k)}  \beta_{\theta}(t) d t .
\end{eqnarray}
From Lie product formula for tensors given by Lemma~\ref{lma: Lie product formula for tensors},  we have 
\begin{eqnarray}\label{eq:117}
 \left\vert \prod\limits_{i=1}^{n} \mathcal{C}_i^{\theta} \right\vert^{\frac{1}{\theta}}
\rightarrow \exp \left(  \sum\limits_{i=1}^{n} \log \mathcal{C}_i  \right). 
\end{eqnarray}
By setting $\theta \rightarrow 0$ in Eqs.~\eqref{eq:114},~\eqref{eq:115} and using Lie product formula given by Eq.~\eqref{eq:117},  we will get Eqs.~\eqref{eq1:thm:Multivaraite Tensor Norm Inequalities} and~\eqref{eq2:thm:Multivaraite Tensor Norm Inequalities}. 
$\hfill \Box$

\section{Tensor Expander Chernoff Bounds Derivation by Majorization}\label{sec:Tensor Expander Chernoff Bounds Derivation by Majorization}



In this section, we will begin with the derivation for the expectation bound of Ky Fan $k$-norm for the product of positive Hermitian tensors in Section~\ref{sec:Expectation Estimation for Product of Tensors}. This bound will play a key role in the next Section~\ref{sec:Tensor Expander Chernoff Bounds} by establishing tensor expander chernoff bounds. 

\subsection{Expectation Estimation for Product of Tensors}\label{sec:Expectation Estimation for Product of Tensors}

The main purpose of this section is to bound the expectation of Ky Fan $k$-norm for the product of positive Hermitian tensors. We extend scalar valued expander Chernoff bound proof in~\cite{healy2008randomness} and matrix valued expander Chernoff bound proof in~\cite{garg2018matrix} to context of tensors and remove the restriction that the summation of all mapped tensors should be zero tensor, i.e., $\sum\limits_{v \in \mathfrak{V}} \mathrm{g}(v) = \mathcal{O}$.






Let $\mathbf{A}$ be the normalized adjacency matrix of the underlying graph $\mathfrak{G}$ and let $\tilde{\mathbf{A}} = \mathbf{A} \otimes \mathcal{I}_{(\mathbb{I}_1^M)^2}$, where the identity tensor $\mathcal{I}_{(\mathbb{I}_1^M)^2}$ has dimensions as $I^2_1 \times \cdots \times I^2_M \times I^2_1 \times \cdots \times I^2_M$. We use $\mathcal{F} \in \mathbb{C}^{\left(n\times I^2_1 \times \cdots \times I^2_M \right)   \times \left(n\times I^2_1 \times \cdots \times I^2_M \right)   }$ to represent block diagonal tensor valued matrix where the $v$-th diagonal block is the tensor
\begin{eqnarray}
\mathcal{T}_v = \exp\left( \frac{t \mathrm{g}(v) (a + \iota b)}{2} \right) \otimes 
 \exp\left( \frac{t \mathrm{g}(v) (a - \iota b)}{2} \right). 
\end{eqnarray}
The tensor $\mathcal{F}$ can also be expressed as 
\begin{eqnarray}\label{eq: block decomposition}
\mathcal{F} &=&\left[
    \begin{array}{cccc}
       \mathcal{T}_{v_1}  &  \mathcal{O}  & \cdots &  \mathcal{O}   \\
        \mathcal{O} &  \mathcal{T}_{v_2}  & \cdots & \mathcal{O}  \\
       \vdots &  \vdots & \ddots & \vdots  \\
       \mathcal{O} &  \mathcal{O} & \cdots & \mathcal{T}_{v_n}  \\
    \end{array}
\right].
\end{eqnarray} 
Then the tensor $\left(\mathcal{F} \star_{M + 1}  \tilde{\mathbf{A}}\right)^{\kappa}$ is a block tensor valued matrix whose $(u, v)$-block is a tensor with dimensions as  $I^2_1 \times \cdots \times I^2_M  \times I^2_1 \times \cdots \times I^2_M$ expressed as :
\begin{eqnarray}
\sum\limits_{v_1, \cdots, v_{\kappa-1} \in \mathfrak{V}} \mathbf{A}_{u,v_1}\left(\prod\limits_{j=1}^{\kappa-2} \mathbf{A}_{v_j, v_{j+1}} \right) 
\mathbf{A}_{v_{\kappa-1}, v} \left(\mathcal{T}_u \star_{2M}  \mathcal{T}_{v_1} \star_{2M} \cdots \star_{2M} \mathcal{T}_{v_{\kappa-1}} \right)
\end{eqnarray}
%

Let $\mathbf{u}_0 \in \mathbb{C}^{n \times \mathbb{I}^2_1 \times \cdots \times \mathbb{I}^2_M}$ be the tensor obtained by $\frac{\mathbf{1}}{\sqrt{n}}\otimes \mbox{\textbf{col}}(\mathcal{I}_{\mathbb{I}_1^M})$, where $\mathbf{1}$ is the all ones vector with size $n$ and $\mbox{\textbf{col}} (\mathcal{I}_{\mathbb{I}_1^M}) \in \mathbb{C}^{I^2_1 \times \cdots \times  I^2_M \times 1 }$ is the column tensor of the identity tensor $\mathcal{I}_{\mathbb{I}_1^M} \in \mathbb{C}^{I_1 \times \cdots \times  I_M \times I_1 \times \cdots \times I_M}$. By applying the following relation:
\begin{eqnarray}\label{eq:kron prod and trace relation}
\bigl<  \mbox{\textbf{col}}(\mathcal{I}_{\mathbb{I}_1^M}),  \mathcal{C} \otimes \mathcal{B}\star_M \mbox{\textbf{col}}(\mathcal{I}_{\mathbb{I}_1^M}) \bigr> = \mathrm{Tr}\left( \mathcal{C} \star_M \mathcal{B}^T \right), 
\end{eqnarray}
where $\mathcal{C}, \mathcal{B} \in \mathbb{C}^{I_1 \times \cdots \times I_M \times I_1 \times \cdots \times I_M}$;
we will have following expectation of $\kappa$ steps transition of Hermitian tensors from the vertex $v_1$ to the vertex $v_{\kappa}$, 
\begin{eqnarray}\label{eq:p17 2nd}
\mathbb{E}\left[\mathrm{Tr} \left(  \prod\limits_{i=1}^{\kappa}\exp \left( \frac{t \mathrm{g}(v_i) (a + \iota b)}{2} \right) \star_{M}  \prod\limits_{i=\kappa}^{1}\exp \left( \frac{t \mathrm{g}(v_i) (a -  \iota b)}{2} \right)\right) \right] = \nonumber \\
=  \mathbb{E}\left[\biggl< \mbox{\textbf{col}}(\mathcal{I}_{\mathbb{I}_1^M}), \prod_{i=1}^\kappa \mathcal{T}_{v_i} \star_M\mbox{\textbf{col}}(\mathcal{I}_{\mathbb{I}_1^M}) \biggr> \right]= \biggl<  \mathbf{u}_0, \left( \mathcal{F}\star_{M+1} \tilde{\mathbf{A}} \right)^{\kappa} \star_{M+1} \mathbf{u}_0  \biggr>.
\end{eqnarray}
If we define $\left( \mathcal{F}\star_{M+1} \tilde{\mathbf{A}} \right)^{\kappa} \star_{M+1} \mathbf{u}_0$ as $\mathbf{u}_{\kappa}$, the goal of this section is to estimate $ \bigl<  \mathbf{u}_0, \mathbf{u}_{\kappa}  \bigr>$. 

The trick is to separate the space of $\mathbf{u}$ as the subspace spanned by the $(\mathbb{I}_1^M)^2$ tensors $\mathbf{1} \otimes e_i$ denoted by $\mathbf{u}^{\parallel}$, where $1 \leq i \leq (\mathbb{I}_1^M)^2$ and $e_i \in \mathbb{C}^{I^2_1 \times \cdots \times I^2_M \times 1 }$ is the column tensor of size $(\mathbb{I}_1^M)^2$ with 1 in position $i$ and 0 elsewhere, and its orthogonal complement space, denoted by $\mathbf{u}^{\perp}$. Following lemma is required to bound how the tensor norm is changed in terms of aforementioned subspace and its orthogonal space after acting by the tensor $\mathcal{F}\star_{2M+1} \tilde{\mathbf{A}}$. 

\begin{lemma}\label{lma:4_4}
Given paramters $\lambda \in (0, 1)$, $a  \geq 0$, $r > 0$, and $t > 0$. Let $\mathfrak{G} = (\mathfrak{V}, \mathfrak{E})$ be a regular $\lambda$-expander graph on the vetices set $\mathfrak{V}$ and $\left\Vert \mathrm{g}(v_i) \right\Vert \leq r$ for all $v_i \in \mathfrak{V}$. Each vertex $v \in \mathfrak{V}$ will be assigned a tensor $\acute{\mathcal{T}}_v$, where $\acute{\mathcal{T}}_v \define \frac{ \mathrm{g}(v) (a + \iota b)}{2} \otimes \mathcal{I}_{\mathbb{I}_1^M}+ 
 \mathcal{I}_{\mathbb{I}_1^M} \otimes \frac{ \mathrm{g}(v) (a - \iota b)}{2} \in \mathbb{C}^{I^2_1 \times \cdots \times I^2_M \times I^2_1 \times \cdots \times  I^2_M}$. Let $\mathcal{F} \in \mathbb{C}^{\left(n\times I^2_1 \times \cdots \times I^2_M \right)   \times \left(n\times I^2_1 \times \cdots \times I^2_M \right)   }$ to represent block diagonal tensor valued matrix where the $v$-th diagonal block is the tensor $\exp(t \acute{\mathcal{T}}_v) = \mathcal{T}_v$. For any tensor $\mathbf{u} \in \mathbb{C}^{n \times \mathbb{I}^2_1 \times \cdots \times \mathbb{I}^2_M}$, we have
\begin{enumerate}
   \item $\left\Vert \left( \mathcal{F} \star_{M+1} \tilde{\mathbf{A}} \star_{M+1} \mathbf{u}^{\parallel} \right)^{\parallel} \right\Vert \leq \gamma_1 \left\Vert  \mathbf{u}^{\parallel} \right\Vert$, where $\gamma_1 = \exp (tr \sqrt{a^2 + b^2}  )$;
   \item $\left\Vert \left( \mathcal{F} \star_{M+1} \tilde{\mathbf{A}} \star_{M+1} \mathbf{u}^{\perp} \right)^{\parallel} \right\Vert \leq \gamma_2 \left\Vert  \mathbf{u}^{\perp} \right\Vert$, where $\gamma_2 = \lambda (\exp(tr \sqrt{a^2 + b^2}) - 1)$;
   \item $\left\Vert \left( \mathcal{F} \star_{M+1} \tilde{\mathbf{A}} \star_{M+1} \mathbf{u}^{\parallel} \right)^{\perp} \right\Vert \leq \gamma_3 \left\Vert  \mathbf{u}^{\parallel} \right\Vert$, where $\gamma_3 = \exp(tr \sqrt{a^2 + b^2})-1$;
   \item $\left\Vert \left( \mathcal{F} \star_{M+1} \tilde{\mathbf{A}} \star_{M+1} \mathbf{u}^{\perp} \right)^{\perp} \right\Vert \leq \gamma_4 \left\Vert  \mathbf{u}^{\perp} \right\Vert$, where $\gamma_4 = \lambda \exp(tr \sqrt{a^2 + b^2})$.
\end{enumerate}
\end{lemma}
\textbf{Proof:}

For Part 1, let $\mathbf{1} \in \mathbb{C}^n$ be all ones vector, and let $\mathbf{u}^{\parallel} =  \mathbf{1} \otimes \mathbf{v}$ for some $\mathbf{v} \in \mathbb{C}^{(\mathbb{I}_1^M)^2}$. Then, we have
\begin{eqnarray}
\left( \mathcal{F} \star_{M+1} \tilde{\mathbf{A}} \star_M \mathbf{u}^{\parallel} \right)^{\parallel} 
= \left( \mathcal{F}  \star_M \mathbf{u}^{\parallel} \right)^{\parallel}
= \mathbf{1} \otimes \left(\frac{1}{n} \sum\limits_{v \in \mathfrak{V}} \exp(t \acute{\mathcal{T}}_v) \star_M \mathbf{v}  \right)
\end{eqnarray}
and we can bound $\left(\frac{1}{n} \sum\limits_{v \in \mathfrak{V}} \exp( t \acute{\mathcal{T}}_v ) \star_M \mathbf{v}  \right)$ further as 
\begin{eqnarray}
\left\Vert \frac{1}{n} \sum\limits_{v \in \mathfrak{V}} \exp(t \acute{\mathcal{T}}_v)  \right\Vert &=& \left\Vert \frac{1}{n}
\sum\limits_{v \in \mathfrak{V}} \sum\limits_{i=0}^{\infty} \frac{t^i \acute{\mathcal{T}}^i_v }{i !}  \right\Vert \nonumber \\
&=& \left\Vert \mathcal{I} + \frac{1}{n}\sum\limits_{v \in \mathfrak{V}} \sum\limits_{i=1}^{\infty} 
\frac{t^i \acute{\mathcal{T}}^i_v }{i !} \right\Vert \nonumber \\
&\leq&1 + \frac{1}{n} \sum\limits_{v \in \mathfrak{V}} \sum\limits_{i=1}^{\infty} 
\frac{t^i \left\Vert \acute{\mathcal{T}}_v \right\Vert^i }{i !} \nonumber \\
&\leq & 1 + \sum\limits_{i=1}^{\infty}\frac{(tr \sqrt{a^2 + b^2})^i }{i !}=\exp (tr \sqrt{a^2 + b^2}  ),
\end{eqnarray}
where the last inequality is due to the fact that $\left\Vert \frac{t \mathrm{g}(v) (a + \iota b)}{2} \otimes \mathcal{I}_{\mathbb{I}_1^M}+ 
 \mathcal{I}_{\mathbb{I}_1^M} \otimes \frac{ t \mathrm{g}(v) (a - \iota b)}{2} \right\Vert \leq 2 tr \times \sqrt{\frac{a^2 + b^2}{4}}  $. 

Then Part 1. of this lemma is established due to 
\begin{eqnarray}
\left\Vert  \left( \mathcal{F} \star_{M+1} \tilde{\mathbf{A}} \star_{M+1} \mathbf{u}^{\parallel} \right)^{\parallel} \right\Vert & =& \sqrt{n} \left\Vert \frac{1}{n} \sum\limits_{v \in \mathfrak{V}}\exp(t \acute{\mathcal{T}}_v) \mathbf{v} \right\Vert \nonumber \\
&\leq & \sqrt{n} \left\Vert \mathbf{v}  \right\Vert \exp (tr \sqrt{a^2 + b^2}  ) =\exp (tr \sqrt{a^2 + b^2}  )  \left\Vert  \mathbf{u}^{\parallel} \right\Vert .
\end{eqnarray}

For Part 2, since $(\tilde{\mathbf{A}} \star_{M+1} \mathbf{u}^{\perp})^{\parallel} = 0$, we have 
\begin{eqnarray}
\left\Vert \left(  \mathcal{F} \star_{M+1} \tilde{\mathbf{A}} \star_{M+1} \mathbf{u}^{\perp} \right)^{\parallel} \right\Vert &=& \left\Vert  (  ( \mathcal{F} - \mathcal{I}) \star_{M+1} \tilde{\mathbf{A}} \star_{M+1} \mathbf{u}^{\perp} )^{\parallel} \right\Vert \nonumber \\
&\leq&   \left\Vert   ( \mathcal{F} - \mathcal{I}) \star_{M+1} \tilde{\mathbf{A}} \star_{M+1} \mathbf{u}^{\perp} \right\Vert \nonumber \\
&\leq&  \max\limits_{v \in \mathfrak{V}} \left\Vert  \exp(t \acute{\mathcal{T}}_v) - \mathcal{I} \right\Vert \cdot  \left\Vert \tilde{\mathbf{A}} \star_{M+1} \mathbf{u}^{\perp}   \right\Vert 
\nonumber \\
&\leq&  \max\limits_{v \in \mathfrak{V}} \left\Vert \sum\limits_{i=1}^{\infty} \frac{t^i \acute{\mathcal{T}}^i_v }{i !} \right\Vert \cdot \left\Vert \tilde{\mathbf{A}} \star_{M+1} \mathbf{u}^{\perp}   \right\Vert  \leq (\exp (tr \sqrt{a^2 + b^2}  ) -1) \lambda \left\Vert  \mathbf{u}^{\perp} \right\Vert ,
\end{eqnarray}
where the last inequality uses that the underlying graph $\mathfrak{G}$ is a $\lambda$-expander  graph, i.e., $\left\Vert \mathbf{A} \mathbf{x}\right\Vert \leq \lambda \cdot \left\Vert \mathbf{x}\right\Vert$. Therefore, Part 2 is also valid. 

For Part 3,  because $(\mathbf{u}^{\parallel})^{\perp} = 0$, we have $(\mathcal{F} \star_{M+1} \tilde{\mathbf{A}} \star_{M+1} \mathbf{u}^{\parallel} )^{\perp} =(\mathcal{F} \star_{M+1} \mathbf{u}^{\parallel} )^{\perp} = ( ( \mathcal{F} - \mathcal{I}  ) \star_{M+1} \mathbf{u}^{\parallel }  )^{\perp} $. Then, we can upper bound as 
\begin{eqnarray}
\left\Vert  ( ( \mathcal{F} - \mathcal{I}  ) \star_{M+1} \mathbf{u}^{\parallel }  )^{\perp} \right\Vert &\leq& \left\Vert  ( \mathcal{F} - \mathcal{I}  ) \star_{M+1} \mathbf{u}^{\parallel }  \right\Vert \nonumber \\
&\leq& \left\Vert  ( \mathcal{F} - \mathcal{I}   \right\Vert \cdot \left\Vert \mathbf{u}^{\parallel }  \right\Vert \nonumber \\
& = & \max\limits_{v \in \mathfrak{V}} \left\Vert  \exp(t \acute{\mathcal{T}}_v)  - \mathcal{I} \right\Vert  \cdot \left\Vert \mathbf{u}^{\parallel }  \right\Vert \nonumber \\
& \leq & \max\limits_{v \in \mathfrak{V}} \left\Vert \sum\limits_{i=1}^{\infty} \frac{t^i \acute{\mathcal{T}}^i_v }{i !} \right\Vert \cdot \left\Vert \mathbf{u}^{\parallel}   \right\Vert  \leq (\exp (tr \sqrt{a^2 + b^2}  ) -1)  \left\Vert  \mathbf{u}^{\parallel} \right\Vert ,
\end{eqnarray}
hence, Part 3 is also proved. 

Finally, for Part 4, we have 
\begin{eqnarray}
\left\Vert (\mathcal{F} \star_{M+1} \tilde{\mathbf{A}} \star_{M+1} \mathbf{u}^{\perp} )^{\perp}\right\Vert &\leq& \left\Vert \mathcal{F} \star_{M+1} \tilde{\mathbf{A}} \star_{M+1} \mathbf{u}^{\perp} \right\Vert \nonumber \\
& \leq &  \left\Vert \mathcal{F} \right\Vert \cdot \left\Vert  \tilde{\mathbf{A}} \star_{M+1} \mathbf{u}^{\perp} \right\Vert \leq  \exp (tr \sqrt{a^2 + b^2}  ) \lambda  \left\Vert  \mathbf{u}^{\perp} \right\Vert ,
\end{eqnarray}
where we use $  \left\Vert \mathcal{F} \right\Vert \leq  \exp (tr \sqrt{a^2 + b^2}  ) $ (shown at previous part) and the underlying graph $\mathfrak{G}$ is a $\lambda$-expander graph. 
$\hfill \Box$

In the following, we will apply Lemma~\ref{lma:4_4} to bound the following term provided by Eq.~\eqref{eq:p17 2nd}
\begin{eqnarray}
\biggl<  \mathbf{u}_0, \left( \mathbf{F}\star_{M+1} \tilde{\mathbf{A}} \right)^{\kappa} \star_{M+1} \mathbf{u}_0  \biggr>
\end{eqnarray}
This bound is formulated by the following Lemma~\ref{lma:4_3}

\begin{lemma}\label{lma:4_3}
Let $\mathfrak{G}$ be a regular $\lambda$-expander graph on the vertex set $\mathfrak{V}$, $g : \mathfrak{V} \rightarrow \mathbb{C}^{I_1 \times \cdots \times I_M \times I_1 \times \cdots \times I_M}$, and let $v_1, \cdots, v_{\kappa}$ be a stationary random walk on $\mathfrak{G}$. If $ t r  \sqrt{a^2 + b^2} < 1$ and $\lambda(  2\exp(tr \sqrt{a^2 + b^2}) - 1 ) \leq 1$, we have:
%
\begin{eqnarray}\label{eq1:lma:4_3}
\mathbb{E}\left[\mathrm{Tr} \left(  \prod\limits_{i=1}^{\kappa}\exp \left( \frac{t \mathrm{g}(v_i) (a + \iota b)}{2} \right) \star_{M}  \prod\limits_{i=\kappa}^{1}\exp \left( \frac{t \mathrm{g}(v_i) (a -  \iota b)}{2} \right)\right) \right] \leq  \nonumber \\
\mathbb{I}_1^M \exp \left[ \kappa\left(  2   t r \sqrt{a^2 + b^2} + \frac{8}{1- \lambda} + \frac{16   t r \sqrt{a^2 + b^2} }{1 - \lambda} \right)\right]. ~~~~~~~~~~~~~~~
\end{eqnarray}

\end{lemma}
\textbf{Proof:}
There are two phases for this proof. The first phase is to bound the evolution of tensor norms $\left\Vert \mathbf{u}_i^{\perp}\right\Vert$ and $\left\Vert \mathbf{u}_i^{\parallel}\right\Vert$, respectively. The second phase is to bound $\gamma_i$ for $1 \leq i \leq 4$ in Lemm~\ref{lma:4_4}. We begin with the derivation for the bound $\left\Vert \mathbf{u}_i^{\perp} \right\Vert $, where $ \mathbf{u}_i$ is the output tensor after acting by the tensor $\mathcal{F} \star_{M+1} \tilde{\mathbf{A}}  $ for $i$ times. It is 
\begin{eqnarray}\label{eq2:lma:4_3}
\left\Vert \mathbf{u}_i^{\perp} \right\Vert & = & \left\Vert (  \mathcal{F} \star_{M+1} \tilde{\mathbf{A}} \star_{M+1} \mathbf{u}_{i-1} )^{\perp} \right\Vert \nonumber \\
& \leq & \left\Vert (  \mathcal{F} \star_{M+1} \tilde{\mathbf{A}} \star_{M+1} \mathbf{u}^{\parallel}_{i-1} )^{\perp} \right\Vert + \left\Vert (  \mathcal{F} \star_{M+1} \tilde{\mathbf{A}} \star_{M+1} \mathbf{u}^{\perp}_{i-1} )^{\perp} \right\Vert \nonumber \\
& \leq_1 & \gamma_3 \left\Vert  \mathbf{u}^{\parallel}_{i-1} \right\Vert +
\gamma_4 \left\Vert  \mathbf{u}^{\perp}_{i-1} \right\Vert \nonumber \\
& \leq_2 & (\gamma_3 + \gamma_3 \gamma_4 +  \gamma_3 \gamma^2_4 + \cdots) \max\limits_{j < i}\left\Vert \mathbf{u}_j^{\parallel} \right\Vert \leq \frac{\gamma_3}{1 - \gamma_4} \max\limits_{j < i}\left\Vert \mathbf{u}_j^{\parallel} \right\Vert,
\end{eqnarray}
where $\leq_1$ is obtained from Lemma~\ref{lma:4_4}, $\leq_2$ is obtained by applying the inequality at $\leq_1$ repeatedly. The next task is to bound  $\left\Vert \mathbf{u}_i^{\parallel} \right\Vert $, we have 
\begin{eqnarray}\label{eq3:lma:4_3}
\left\Vert \mathbf{u}_i^{\parallel} \right\Vert & = & \left\Vert (  \mathcal{F} \star_{M+1} \tilde{\mathbf{A}} \star_{M+1} \mathbf{u}_{i-1} )^{\parallel} \right\Vert \nonumber \\
& \leq & \left\Vert (  \mathcal{F} \star_{M+1} \tilde{\mathbf{A}} \star_{M+1} \mathbf{u}^{\parallel}_{i-1} )^{\parallel} \right\Vert + \left\Vert (  \mathcal{F} \star_{M+1} \tilde{\mathbf{A}} \star_{M+1} \mathbf{u}^{\perp}_{i-1} )^{\parallel} \right\Vert \nonumber \\
& \leq_1 & \gamma_1 \left\Vert \mathbf{u}^{\parallel}_{i-1} \right\Vert +
\gamma_2 \left\Vert  \mathbf{u}^{\perp}_{i-1} \right\Vert \nonumber \\
& \leq_2 & \left( \gamma_1 + \frac{\gamma_2 \gamma_3}{ 1 - \gamma_4} \right) \max\limits_{j < i}\left\Vert \mathbf{u}_j^{\parallel} \right\Vert,
\end{eqnarray}
where $\leq_1$ is obtained from Lemma~\ref{lma:4_4}, $\leq_2$ is obtained from Eq.~\eqref{eq2:lma:4_3}. From  Eqs~\eqref{eq:p17 2nd},~\eqref{eq2:lma:4_3} and ~\eqref{eq3:lma:4_3}, we have 
\begin{eqnarray}\label{eq4:lma:4_3}
\mathbb{E}\left[\mathrm{Tr} \left(  \prod\limits_{i=1}^{\kappa}\exp \left( \frac{t \mathrm{g}(v_i) (a + \iota b)}{2} \right) \star_{M}  \prod\limits_{i=\kappa}^{1}\exp \left( \frac{t \mathrm{g}(v_i) (a -  \iota b)}{2} \right)\right) \right]   \nonumber \\
= \langle \mathbf{u}_0, \mathbf{u}_\kappa \rangle 
=  \langle \mathbf{u}_0, \mathbf{u}^{\parallel}_\kappa \rangle  \leq  \left\Vert \mathbf{z}_0 \right\Vert \cdot \left\Vert \mathbf{z}^{\parallel}_\kappa \right\Vert = \sqrt{\mathbb{I}_1^M} \cdot \left\Vert \mathbf{z}^{\parallel}_\kappa \right\Vert  ~~~~~~~~~~~~~~~ \nonumber \\
\leq   \sqrt{\mathbb{I}_1^M} \left( \gamma_1 + \frac{\gamma_2 \gamma_3}{ 1 - \gamma_4} \right)^{\kappa} \cdot \left\Vert \mathbf{z}^{\parallel}_0 \right\Vert 
\leq  \mathbb{I}_1^M \left( \gamma_1 + \frac{\gamma_2 \gamma_3}{ 1 - \gamma_4} 
\right)^{\kappa}. ~~~~~~~~~~
\end{eqnarray}

The second phase of this proof requires us to bound following four terms: $\gamma_i$ for $1 \leq i \leq 4$. Since $tr \sqrt{a^2 + b^2} < 1$, we can bound $\gamma_1$ as following:
\begin{eqnarray}\label{eq:gamma 1}
\gamma_1 &=& \exp (tr \sqrt{a^2 + b^2}) \leq 1 + 2  t r\sqrt{a^2 + b^2};
\end{eqnarray}
\begin{eqnarray}\label{eq:gamma 2}
\gamma_2 &=& \lambda (\exp (tr \sqrt{a^2 + b^2}) - 1) \leq 2 \lambda   t r\sqrt{a^2 + b^2};
\end{eqnarray}
\begin{eqnarray}\label{eq:gamma 3}
\gamma_3 &=& \exp (tr \sqrt{a^2 + b^2}) - 1 \leq 2   t r \sqrt{a^2 + b^2};
\end{eqnarray}
and the condition $\lambda(  2\exp(tr \sqrt{a^2 + b^2}) - 1 ) \leq 1$, we have 
\begin{eqnarray}\label{eq:gamma 4}
1 - \gamma_4 &  = &1 -  \lambda \exp (t r \sqrt{a^2 + b^2} )  \geq \frac{1 -  \lambda}{2}.
\end{eqnarray}
By applying Eqs.~\eqref{eq:gamma 1},~\eqref{eq:gamma 2},~\eqref{eq:gamma 3} and~\eqref{eq:gamma 4} to the upper bound in Eq.~\eqref{eq4:lma:4_3}, we also have 
\begin{eqnarray}
\mathbb{I}_1^M \left( \gamma_1 + \frac{\gamma_2 \gamma_3}{ 1 - \gamma_4} 
\right)^{\kappa} &\leq&  \mathbb{I}_1^M \left[ 1 + 2 (  t r \sqrt{a^2 + b^2} ) + \frac{8 \lambda t^2r^2(a^2 + b^2) }{    1-  \lambda       } \right]^\kappa \nonumber \\
& \leq & \mathbb{I}_1^M \left[ \left(1 + 2  t r \sqrt{a^2 + b^2} \right)  \left(1 + \frac{8}{1 - \lambda} \right) \right]^{\kappa} \nonumber \\
& \leq & \mathbb{I}_1^M \exp \left[ \kappa\left(  2   t r \sqrt{a^2 + b^2} + \frac{8}{1- \lambda} + \frac{16   t r \sqrt{a^2 + b^2} }{1 - \lambda} \right)\right]
\end{eqnarray}
This lemma is proved.
$\hfill \Box$

\subsection{Tensor Expander Chernoff Bounds}\label{sec:Tensor Expander Chernoff Bounds}

We begin with a lemma about a Ky Fan $k$-norm inequality for the sum of tensors.

\begin{lemma}\label{lma:Ky Fan Inequalities for the sum of tensosrs}
Let $\mathcal{C}_i \in \mathbb{C}^{I_1 \times \cdots \times I_N \times I_1 \times \cdots \times I_N}$ with Hermitian rank $r$, then we have 
\begin{eqnarray}\label{eq1:lma:Ky Fan Inequalities for the sum of tensosrs}
\left\Vert \left\vert  \sum\limits_{i=1}^{m} \mathcal{C}_i \right\vert^s \right\Vert_{(k)}
\leq  m^{s -1} \sum\limits_{i=1}^{m}  \left\Vert \left\vert \mathcal{C}_i \right\vert^{s} \right\Vert_{(k)}     
\end{eqnarray}
where $s \geq 1$ and $k \in \{1,2,\cdots, r \}$. 
\end{lemma}
\textbf{Proof:}
Since we have 
\begin{eqnarray}
\left\Vert \left\vert  \sum\limits_{i=1}^{m} \mathcal{C}_i \right\vert^s \right\Vert_{(k)}
 = \sum\limits_{j=1}^{k} \lambda_j \left( \left\vert  \sum\limits_{i=1}^{m} \mathcal{C}_i \right\vert^s \right) =  \sum\limits_{j=1}^{k} \lambda^s_j \left( \left\vert  \sum\limits_{i=1}^{m} \mathcal{C}_i \right\vert \right) = \sum\limits_{j=1}^{k} \sigma^s_j \left( \sum\limits_{i=1}^{m} \mathcal{C}_i \right). 
\end{eqnarray}
where we have orders for eigenvalues as $\lambda_1 \geq \lambda_2 \geq \cdots $, and singular values as  $\sigma_1 \geq \sigma_2 \geq \cdots $. 

From Theorem G.1.d. in~\cite{MR2759813} and Theorem 5.2 in~\cite{ni2019hermitian}, we have Ky Fan singular value majorization inequalities:%
\begin{eqnarray}
\sum\limits_{j = 1}^{k} \sigma_j ( \sum\limits_{i=1}^m \mathcal{C}_i ) \leq 
\sum\limits_{j = 1}^{k} \left( \sum\limits_{i=1}^m \sigma_j (\mathcal{C}_i ) \right),
\end{eqnarray}
where $k \in \{1,2,\cdots, s \}$. Then, we have 
\begin{eqnarray}
\sum\limits_{j = 1}^{k} \sigma^s_j ( \sum\limits_{i=1}^m \mathcal{C}_i ) &\leq &
\sum\limits_{j = 1}^{k} \left( \sum\limits_{i=1}^m \sigma_j (\mathcal{C}_i ) \right)^s
\leq m^{s-1} \sum\limits_{j = 1}^{k} \left( \sum\limits_{i=1}^m \sigma^s_j (\mathcal{C}_i ) \right) \nonumber \\
& = & 
m^{s-1} \sum\limits_{j = 1}^{k} \left( \sum\limits_{i=1}^m \sigma^s_j ( \left\vert \mathcal{C}_i \right\vert ) \right) = m^{s-1} \sum\limits_{j = 1}^{k} \left( \sum\limits_{i=1}^m \sigma_j ( \left\vert \mathcal{C}_i \right\vert^s ) \right) \nonumber \\
& = & m^{s -1} \sum\limits_{i=1}^{m}  \left\Vert \left\vert \mathcal{C}_i \right\vert^{s} \right\Vert_{(k)}    
\end{eqnarray}
$\hfill \Box$

We are ready to present our main theorem about the tensor expander bound for Ky Fan $k$-norm.
\begin{theorem}\label{thm:tensor expander}
Let $\mathfrak{G} = (\mathfrak{V}, \mathfrak{E})$ be a regular undirected graph whose transition matrix has second eigenvalue $\lambda$, and let $g: \mathfrak{V} \rightarrow \in \mathbb{C}^{I_1 \times \cdots \times I_M \times I_1 \times \cdots \times I_M}$ be a function. We assume following: 
\begin{enumerate}
\item For each $v \in \mathfrak{V}$, $g(v)$ is a Hermitian tensor;
\item $\left\Vert g(v) \right\Vert \leq r$;
\item A nonnegative coefficients polynomial raised by the power $s \geq 1$ as $f: x \rightarrow (a_0 + a_1x +a_2 x^2 + \cdots +a_n x^n)^s$ satisfying $f \left(\exp \left( t   \sum\limits_{j=1}^{\kappa} g(v_j) \right) \right) \succeq \exp \left( t f \left(  \sum\limits_{j=1}^{\kappa} g(v_j) \right) \right) $ almost surely;
\item  For $\tau \in [\infty, \infty]$, we have constants $C$ and $\sigma$ such that $ \beta_0(\tau)  \leq \frac{C \exp( \frac{-\tau^2}{2 \sigma^2} ) }{\sigma \sqrt{2 \pi}}$. 
\end{enumerate}
Then, we have 
\begin{eqnarray}\label{eq0:thm:tensor expander}
\mathrm{Pr}\left( \left\Vert f \left(  \sum\limits_{j=1}^{\kappa} g(v_j) \right)  \right\Vert_{(k)} \geq\vartheta \right) \leq   \min\limits_{t > 0 } \left[ (n+1)^{(s-1)} e^{-\vartheta t} \left(a_0 k  +C \left( k + \sqrt{\frac{\mathbb{I}_1^M - k }{k}}\right)\cdot \right. \right. ~~~~~  \nonumber \\
\left. \left. \sum\limits_{l=1}^n a_l \exp( 8 \kappa \overline{\lambda} + 2  (\kappa +8 \overline{\lambda}) lsr t + 2(\sigma (\kappa +8 \overline{\lambda}) lsr )^2 t^2  )  \right)\right],
\end{eqnarray}
where $\overline{\lambda} = 1 - \lambda$. 
\end{theorem}
\textbf{Proof:}
Let $t > 0$ be a paramter to be chosen later, then we have 
\begin{eqnarray}\label{eq1:thm:tensor expander}
\mathrm{Pr}\left( \left\Vert f \left(  \sum\limits_{j=1}^{\kappa} g(v_j) \right)  \right\Vert_{(k)} \geq \vartheta \right) &=&\mathrm{Pr}\left( \exp\left( \left\Vert t f \left(  \sum\limits_{j=1}^{\kappa} g(v_j) \right)  \right\Vert_{(k)} \right) \geq \exp\left(\vartheta t \right) \right) \nonumber \\
&=_1& \mathrm{Pr}\left(  \left\Vert \exp\left( t f \left(  \sum\limits_{j=1}^{\kappa} g(v_j) \right)  \right)  \right\Vert_{(k)} \geq \exp\left(\vartheta t \right) \right) \nonumber \\
&\leq_2 &   \exp\left( - \vartheta t \right) \mathbb{E} \left(  \left\Vert \exp\left( t f \left(  \sum\limits_{j=1}^{\kappa} g(v_j) \right)  \right)  \right\Vert_{(k)}  \right) \nonumber \\
&\leq_3 &  \exp\left( - \vartheta t \right) \mathbb{E} \left(  \left\Vert f \left( \exp\left( t   \sum\limits_{j=1}^{\kappa} g(v_j) \right)  \right)  \right\Vert_{(k)}  \right) ,
\end{eqnarray}
where equality $=_1$ comes from spectral mapping theorem, inequality $\leq_2$ is obtained from Markov inequality, and the last inequality $\leq_3$ is based on our function $f$ assumption (third assumption).

From Eq.~\eqref{eq2:thm:Multivaraite Tensor Norm Inequalities} in Theorem~\ref{thm:Multivaraite Tensor Norm Inequalities}, we can further bound the expectation term in Eq.~\eqref{eq1:thm:tensor expander} as 
\begin{eqnarray}\label{eq2:thm:tensor expander}
\mathbb{E} \left(  \left\Vert f \left( \exp\left( t   \sum\limits_{j=1}^{\kappa} g(v_j) \right)  \right)  \right\Vert_{(k)}  \right)  ~~~~~~~~~~~~~~~~~~~~~~~~~~~~~~~~~~~~~~~~~~~~~~~~~~~~~~~~~ \nonumber \\
 \leq \mathbb{E} \left( \int\limits_{- \infty}^{\infty} \left\Vert f \left(   \left\vert \prod\limits_{j=1}^{\kappa} \exp\left( t g(v_j) (1 + \iota \tau   ) \right) \right\vert  \right) \right\Vert_{(k)}\beta_0( \tau) d   \tau   \right) ~~~~~~~~~~~~~~~~~~~~~~~~~~~ \nonumber \\
=_1  \mathbb{E} \left( \int\limits_{- \infty}^{\infty} \left\Vert \left(
\sum\limits_{l=0}^{n}a_l   \left\vert \prod\limits_{j=1}^{\kappa} \exp\left( t g(v_j) (1 + \iota \tau   ) \right) \right\vert^l 
 \right)^s
 \right\Vert_{(k)}  \beta_0( \tau) d \tau  \right)  ~~~~~~~~~~~~~~~~~~ \nonumber \\
\leq_2 (n+1)^{(s-1)} \mathbb{E} \left( \int\limits_{- \infty}^{\infty} \sum\limits_{l=0}^{n}a_l  \left\Vert \left(
  \left\vert \prod\limits_{j=1}^{\kappa} \exp\left( t g(v_j) (1 + \iota \tau   ) \right) \right\vert^l 
 \right)^s
 \right\Vert_{(k)} \beta_0( \tau) d   \tau   \right)  ~~  \nonumber \\
= (n+1)^{(s-1)} \cdot ~~~~~~~~~~~~~~~~~~~~~~~~~~~~~~~~~~~~~~~~~~~~~~~~~~~~~~~~~~~~~~~~~~~~~~~~~~~~~~ \nonumber \\
\left( \int\limits_{- \infty}^{\infty} \sum\limits_{l=0}^{n}a_l  \mathbb{E} \left( \left\Vert \left(
  \left\vert \prod\limits_{j=1}^{\kappa} \exp\left( t g(v_j) (1 + \iota \tau   ) \right) \right\vert^l 
 \right)^s
 \right\Vert_{(k)} \right)  \beta_0( \tau) d   \tau   \right) 
\end{eqnarray}
where equality $=_1$ comes from the function $f$ definition, inequality $\leq_2$ is based on Lemma~\ref{lma:Ky Fan Inequalities for the sum of tensosrs}. Each summand for $l \ge 1$ in Eq.~\eqref{eq2:thm:tensor expander} can further be bounded as 
\begin{eqnarray}\label{eq3:thm:tensor expander}
 \mathbb{E} \left( \left\Vert \left(
  \left\vert \prod\limits_{j=1}^{\kappa} \exp\left( t g(v_j) (1 + \iota \tau   ) \right) \right\vert^l 
 \right)^s
 \right\Vert_{(k)} \right)  = \mathbb{E} \left( \left\Vert 
  \left\vert \prod\limits_{j=1}^{\kappa} \exp\left( t g(v_j) (1 + \iota \tau   ) \right) \right\vert^{ls} 
  \right\Vert_{(k)} \right) \nonumber \\
\leq_1   \frac{k}{\mathbb{I}_1^M} \mathrm{Tr}\left(   \left\vert \prod\limits_{j=1}^{\kappa} \exp\left( t g(v_j) (1 + \iota \tau   ) \right) \right\vert^{ls} \right) + ~~~~~~~~~~~~~~~~~~~~~~~~~~~~~~~~~~~~~~~~~~~ \nonumber \\
  \sqrt{\frac{\mathbb{I}_1^M - k}{k \mathbb{I}_1^M }   \mathrm{Tr} \left(  \left\vert \prod\limits_{j=1}^{\kappa} \exp\left( t g(v_j) (1 + \iota \tau   ) \right) \right\vert^{2ls} \right)} ~~~~~~~~~~~~~~~~~~~~~~~~~~~~~~~~~~~~ \nonumber \\
\leq_2 k 
\exp \left[ \kappa\left(  2   lst r \sqrt{1 + \tau^2} + \frac{8}{1- \lambda} + \frac{16   tls r \sqrt{1 + \tau^2} }{1 - \lambda} \right)\right] + ~~~~~~~~~~~~~~~~~~~~~ \nonumber \\
\left( \frac{\mathbb{I}_1^M - k}{k} \exp \left[ \kappa\left(  4 lst r \sqrt{1 + \tau^2} + \frac{8}{1- \lambda} + \frac{32  l s t r \sqrt{1 + \tau^2} }{1 - \lambda} \right)\right] \right)^{1/2} ~~~~~~ \nonumber \\
\leq_3 \left( k + \sqrt{\frac{\mathbb{I}_1^M - k }{k}}\right)\cdot \exp\left[ \kappa\left(  2   lst r (1 +\tau) + \frac{8}{1- \lambda} + \frac{16   tls r (1 +\tau)  }{1 - \lambda} \right)\right] ~~~~
\end{eqnarray}
where $\leq_1$ comes from Theorem 5 in~\cite{merikoski1997bounds}, and $\leq_2$ comes from Lemma~\ref{lma:4_3}, and the last inequality $\leq_3$ is obtaiend by bounding $\sqrt{1 + \tau^2}$ as $1 + \tau$.

From Eqs.~\eqref{eq1:thm:tensor expander},~\eqref{eq2:thm:tensor expander}, and~\eqref{eq3:thm:tensor expander}, we have 
\begin{eqnarray}\label{eq4:thm:tensor expander}
\mathrm{Pr}\left( \left\Vert f \left(  \sum\limits_{j=1}^{\kappa} g(v_j) \right)  \right\Vert_{(k)} \geq \vartheta \right) \leq \min\limits_{t > 0   } \left[ (n+1)^{(s-1)} e^{-\vartheta t} \left(a_0 k  + \left( k + \sqrt{\frac{\mathbb{I}_1^M - k }{k}}\right)\cdot \right. \right.  \nonumber \\
\left. \left. \sum\limits_{l=1}^n a_l \int\limits_\infty^{\infty} 
\exp\left[ \kappa\left(  2   lst r (1 +\tau) + \frac{8}{1- \lambda} + \frac{16   tls r (1 +\tau)  }{1 - \lambda} \right)\right] \beta_0(\tau) d \tau \right)\right] ~~~~~~~ \nonumber \\
\leq_1  \min\limits_{t > 0   } \left[ (n+1)^{(s-1)} e^{-\vartheta t} \left(a_0 k  + \left( k + \sqrt{\frac{\mathbb{I}_1^M - k }{k}}\right)\cdot \right. \right.  ~~~~~~~~~~~~~~~~~~~~~~~~~~~~~~~~~~~~~~~ \nonumber \\
\left. \left. \sum\limits_{l=1}^n a_l \int\limits_\infty^{\infty} 
\exp\left[ \kappa\left(  2   lst r (1 +\tau) + \frac{8}{1- \lambda} + \frac{16   tls r (1 +\tau)  }{1 - \lambda} \right)\right]  \frac{C \exp( \frac{-\tau^2}{2 \sigma^2} ) }{\sigma \sqrt{2 \pi}} d \tau  \right)\right] \nonumber \\
 =  \min\limits_{t > 0 } \left[ (n+1)^{(s-1)} e^{-\vartheta t} \left(a_0 k  +C \left( k + \sqrt{\frac{\mathbb{I}_1^M - k }{k}}\right)\cdot \right. \right. ~~~~~~~~~~~~~~~~~~~~~~~~~~~~~~~~~~~~~ \nonumber \\
\left. \left. \sum\limits_{l=1}^n a_l \exp( 8 \kappa \overline{\lambda} + 2  (\kappa +8 \overline{\lambda}) lsr t + 2(\sigma (\kappa +8 \overline{\lambda}) lsr )^2 t^2  )  \right)\right],~~~~~~~~~~~~~~~~~~~
\end{eqnarray}
where inequality $\leq_1$ is obtained by the distribution bound for $ \beta_0(\tau)$ via another distribution function $\frac{C \exp( \frac{-\tau^2}{2 \sigma^2} ) }{\sigma \sqrt{2 \pi}}$, and the last equality comes from Gaussian integral with respect to the variable $\tau$ by setting $1 -\lambda$ as $\overline{\lambda}$. 
$\hfill \Box$

Following corollary is about a tensor expander bound with identity function $f$. 

\begin{corollary}\label{cor:tensor expander one variable}
If we consider the special case of Theorem~\ref{thm:tensor expander} by assuming that the function $f: x \rightarrow x$ is an identiy map, then we have 
\begin{eqnarray}
\mathrm{Pr}\left( \left\Vert  \sum\limits_{j=1}^{\kappa} g(v_j)  \right\Vert_{(k)} \geq  \vartheta   \right) &\leq &  C \left( k + \sqrt{\frac{\mathbb{I}_1^M - k }{k}} \right) \cdot \nonumber \\
&  & \exp \left( - \frac{ \vartheta ^2}{8 \sigma^2 r^2} + \frac{ \vartheta }{2 \sigma^2 r^2} - \frac{1}{2\sigma^2} + 8\kappa \overline{\lambda} \right).
\end{eqnarray}
\end{corollary}
\textbf{Proof:}
From Theorem~\ref{thm:tensor expander}, since the exponent is a quadratic function of $t$, the minimum of this quadratic function is achieved by selecting $t$ as
\begin{eqnarray}\label{eq1:cor:tensor expander one variable}
t = \frac{\vartheta - 2 (\kappa + 8 \overline{\lambda}) r}{4 \sigma^2 r^2 (\kappa + 8 \overline{\lambda} )^2},
\end{eqnarray}
then, we have the desired bound after some algebra by applying Eq.~\eqref{eq1:cor:tensor expander one variable} in Eq.~\eqref{eq0:thm:tensor expander} and setting $l=s=1$, all $a_i=0$ for $1 \leq i \leq n$ except $a_1 = 1$. 
$\hfill \Box$

\section{Conclusions}\label{sec:Conclusions}

In this work, we first build tensor norm inequalities based on the concept of log-majorization, and apply these new tensor norm inequalities to derive the \emph{tensor expander Chernoff bounds} which generalize the matrix expander Chernoff bound by adopting more general norm for tensors, Ky Fan norm, and general convex function, instead of identity function, of random tensors summation. 

There are several future directions that can be explored based on the current work. The first is to consider other types of tensor expander Chernoff bounds under other non-independent assumptions among random tensors. The other direction is to characterize random behaviors of other tensor related quantities besides norms or eigenvalues, for example, what is the Hermitian tensor rank behavior for the summation of random tensors.

\section*{Acknowledgments}
The helpful comments of the referees are gratefully acknowledged.

\bibliographystyle{siamplain}
\bibliography{TensorExpanderChernoff_Bib}
\end{document}